\documentclass[reqno,twoside,11pt,english]{amsart}
\usepackage{amsmath,amsfonts,amssymb,amsthm,epsfig}
\newcommand{\RNum}[1]{\uppercase\expandafter{\romannumeral #1\relax}}
\usepackage{comment}
\usepackage{color}
\usepackage[final,colorlinks,linkcolor=blue,anchorcolor=red,citecolor=blue]{hyperref}
\usepackage{graphicx}
\usepackage{titletoc,epsf}
\usepackage{amsmath,amsfonts,latexsym,amsthm,amsxtra,amssymb,bbm}
\allowdisplaybreaks
\usepackage{graphicx,float}

\newcommand{\Si}{\Sigma}
\newcommand{\lam}{\lambda}
\newcommand{\La}{\Lambda}

\newcommand{\ba}{\backslash}

\newcommand{\la}{\lambda}
\newcommand{\de}{\delta}
\newcommand{\De}{\Delta}
\newcommand{\na}{\nabla}
\newcommand{\Om}{\Omega}
\newcommand{\om}{\omega}
\newcommand{\oo}{\infty}

\newcommand{\pa}{\partial}

\newcommand{\non}{\nonumber}

\newcommand{\al}{\alpha}
\newcommand{\dist}{\text{\rm dist}}
\newcommand{\medint}{-\kern -,375cm\int}         
\newcommand{\medintinrigo}{-\kern -,315cm\int}
\newcommand{\wto}{\rightharpoonup}                

  \newcommand{\wq}{\infty}
 \newcommand{\cal}{\mathcal}

\newcommand{\3}{{\mathcal C}}

\newcommand{\HH}{{\mathcal H}}

\newcommand{\R}{{\mathbb R}}

\newcommand{\N}{{\mathbb N}}

\newcommand{\ep}{\varepsilon}
\newcommand{\su}{\subset}

\newcommand{\ri}{\rightarrow}
\newcommand{\rii}{\rightharpoonup}
\def\supp{\textup{supp}}

\def\loc{\text{\rm loc}}

\newcommand{\divg}{\text{\rm div}}
\newcommand{\reg}{\text{\rm reg}}
\newcommand{\sing}{\text{\rm sing}}
\newcommand{\spt}{\text{\rm spt}}

\newcommand{\Span}{\text{\rm span}}

\numberwithin{equation}{section}
\newtheorem{theorem}{Theorem}[section]

\newtheorem{corollary}[theorem]{Corollary}

\newtheorem{definition}[theorem]{Definition}

\newtheorem{lemma}[theorem]{Lemma}

\newtheorem{proposition}[theorem]{Proposition}
\newtheorem{remark}[theorem]{Remark}
\newtheorem*{remark*}{Remark}

\theoremstyle{definition}

\makeatother

\begin{document}

\title[Biharmonic maps]{Optimal higher regularity for biharmonic maps via quantitative stratification}

 \author[C.-Y. Guo, G.-C. Jiang, C.-L. Xiang and G.-F. Zheng ]{Chang-Yu Guo, Gui-Chun Jiang, Chang-Lin Xiang$^*$ and Gao-Feng Zheng}

\address[Chang-Yu Guo]{Research Center for Mathematics and Interdisciplinary Sciences, Shandong University 266237,  Qingdao, P. R. China and and  Frontiers Science Center for Nonlinear Expectations, Ministry of Education, P. R. China}
\email{changyu.guo@sdu.edu.cn}

\address[Gui-Chun Jiang]{School of Mathematics and Statistics, Hubei Normal University, Huangshi 435002, P. R. China}
\email{gcjiang@hbnu.edu.cn}

\address[Chang-Lin Xiang]{Three Gorges Mathematical Research Center, China Three Gorges University,  443002, Yichang,  P. R. China}
\email{changlin.xiang@ctgu.edu.cn}

\address[Gao-Feng Zheng]{School of Mathematics and Statistics, and Key Lab NAA--MOE, Central China Normal University, Wuhan 430079,  P. R.  China}
\email{gfzheng@ccnu.edu.cn}

\thanks{$^*$Corresponding author: Chang-Lin Xiang}
\thanks{C.-Y. Guo is supported by the Young Scientist Program of the Ministry of Science and Technology
	of China (No.~2021YFA1002200), the National Natural Science Foundation of China (No.~12101362), the Taishan Scholar Project and the Natural Science Foundation of Shandong Province (No.~ZR2022YQ01, ZR2021QA003). The corresponding author C.-L. Xiang is financially supported by the National Natural Science Foundation of China (No.~12271296), and partly supported by the Open Research Fund of Key Laboratory of Nonlinear Analysis \& Applications  (Central China Normal University), Ministry of Education, P. R. China. G.-F. Zheng is supported by the National Natural Science Foundation of China (No.~12271195 and No.~12171180). }

\date{}

\begin{abstract}
This little note is devoted to refining the almost optimal regularity results of Breiner and Lamm \cite{Breiner-Lamm-2015} on minimizing and stationary biharmonic maps  via the powerful quantitative stratification method introduced by Cheeger and Naber \cite{Cheeger-Naber-2013} and further developed by Naber and Valtorta \cite{Naber-V-2017,Naber-V-2018}  for harmonic maps. In particular, we obtain an optimal regularity results for minimizing biharmonic maps. 
\end{abstract}



\maketitle

{\small
\smallskip
\newline
	\tableofcontents}
\bigskip

\section{Introduction and main results}
Let $\Omega\subset\mathbb{R}^m$ be a bounded smooth domain with $m\geq 4$, and let $N\subset \mathbb{R}^L$ be an isometrically embedded compact submanifold  without boundary.  In this note we are concerned with local analysis of (extrinsic) biharmonic maps, which are defined as critical points of  the bienergy functional
$$E(u):=\frac{1}{2}\int_\Omega|\De u|^2dx, \qquad u\in H^2(\Om, N),$$ where
\[H^2(\Om, N)=\Big\{u\in H^2(\Om, \mathbb{R}^L):  u(x)\in N \text{ for almost every }x\in N\Big\}.\]
More precisely,
\begin{definition}\label{uy}
	$u\in H^{2}(\Omega,N)$ is called a weakly biharmonic map if
	\begin{eqnarray}\label{1.1}
	\frac{d}{dt}\Big|_{t=0}E(\pi_N(u+t\xi))=0,\qquad \text{for all }\,\xi\in C_0^\infty(\Omega,\mathbb{R}^L),
	\end{eqnarray}
	where $\pi_N$ is the nearest point projection with respect to $N$.
	
	If a weakly harmonic map $u\colon \Omega\to N$ also satisfies
	\begin{eqnarray}\label{1.1}
	\frac{d}{dt}\Big|_{t=0}E\Big(u(x+t\psi(x))\Big)=0,\qquad \text{for all }\,\psi\in C_0^\infty(\Omega,\mathbb{R}^m),
	\end{eqnarray}
	then it is called a stationary biharmmonic map.
	
	Finally, $u\colon \Omega\to N$ is  a minimizing biharmonic map if
	\begin{eqnarray}\label{1.7}
	E(u)\leq E(v)
	\end{eqnarray}
	for all $v\in H^2(\Omega,N)$ with $v=u$ on $\Omega\backslash K$ for some compact $K\subset\Omega$.
\end{definition}
It should be noticed that the bienergy depends on the embedding $N\hookrightarrow \R^L$. For this reason, biharmonic maps considered in this paper are also called \emph{extrinsic} biharmonic maps in literature (see for instance \cite{Wang-2004-CPAM,Lamm-Riviere-2008}). All the regularity theory studied in this paper  can be  similarly considered for the so-called \emph{intrinsic} biharmonic maps, see e.g.  \cite{ Lamm-Riviere-2008, Moser-2008-CPDE, Scheven-2009-Poincare,Wang-2004-CPAM}.

Motivated by problems from conformal geometry, S.-Y. Chang, L. Wang and P. Yang \cite{Chang-W-Y-1999} introduced the class of biharmonic maps and proved that every weakly biharmonic map $u\in H^1(\Om,S^k)$ is smooth if the dimension $m=4$, and every stationary biharmonic map satisfies the Hausdorff dimension bound $\dim_{\HH}\sing(u)\leq m-4$ if $m\geq 5$, where the singular set $\sing(u)$ of $u$ is defined as
\[
\sing(u):=\{x\in \Omega: u \text{ is not continuous at }x\}.
\]
Later,  C.Y. Wang \cite{Wang-2004-CVPDE} and P. Strzelecki \cite{Strzelecki-2003-CV} independently provided alternative proofs of the above regularity result. A similar regularity result for biharmonic maps into general Riemannian manifold $N$ was achieved soon by C.Y. Wang in the pioneer works \cite{Wang-2004-CPAM,Wang-2004-MathZ}. This can be seen as a biharmonic analogue of the famous regularity theories of  \cite{Bethuel-1993,Evans-1991,Schoen-Uhlenbeck-1982} for harmonic maps. All the above regularity results on biharmonic maps relies on the idea of H\'elein \cite{Helein-1990}: exploiting the hidden symmetry of the target manifold. Different approaches for these regularity results were provided later by T. Lamm and T. Rivi\`ere \cite{Lamm-Riviere-2008}, M. Struwe \cite{Struwe-2008} and \cite{Guo-Xiang-2020-JLMS,Guo-Xiang-Zheng-2021-CVPDE,Guo-Wang-Xiang-2023-CVPDE}, relying on the fundamental works of T. Rivi\`ere \cite{Riviere-2007-Invent} and Rivi\`ere-Struwe \cite{Rivere-Struwe-2008}. In another interesting work, C. Scheven \cite{Scheven-2008-ACV} applied the analysis of defect measures, which was first developed by F.H. Lin for harmonic maps in his seminal work \cite{Lin-1999-Annals}, to obtain an improved dimensional estimate $\dim_{\HH}\sing(u)\leq m-5$ for minimizing biharmonic maps. Furthermore, if the target manifolds $N$ does not admit nonconstant harmonic sphere $u\in W^{2,2}(S^k,N)$ for $4\leq k\leq k_0$, then one has a better bound $\dim_{\HH}\sing(u)\leq m-k_0-2$. For other aspects of biharmonic maps, such as heat flow and bubbling analysis, see for instance \cite{Angelsberg-2006,Chen-Zhu-2023,Lamm-2005,Liu-Yin-2016,Wang-Zheng-2012-JFA}. We would like to point out that, in contrast to harmonic map \cite{Lin-1999-Annals,Simon-book-1996}, the classical stratification theory for biharmonic maps remains merely partially known \cite{Scheven-2008-ACV}.  

In a recent major breakthrough, J. Cheeger and A. Naber \cite{Cheeger-Naber-2013} introduced a powerful new technique, known as \emph{quantitative stratification} in literature, and sucessfully applied it in a fundamental work \cite{Cheeger-Naber-2013-CPAM} to obtain almost optimal regularity for minimizing harmonic maps and minimal currents; see also \cite{Cheeger-H-N-2013,Cheeger-H-N-2015,Cheeger-N-V-2015,Cheeger-J-N-2021,Naber-V-2020-JEMS} for other successful applications of this technique. Soon after Cheeger-Naber's seminal work \cite{Cheeger-Naber-2013-CPAM} on minimizing harmonic maps, C. Breiner and T. Lamm \cite{Breiner-Lamm-2015} extended the quantitative stratificaiton technique to the class of minimizing biharmonic maps. To record their main results, we need to introduce a few concepts. 

For a measurable map $u\colon B_3(0)\subset \R^m\to N$, $r$, $\epsilon>0$ and $k>0$, we may define the \emph{$k$-th quantitative stratum} $S_{\epsilon,r}^k(u)\subset B_3(0)$ (see Definition \ref{ddd} below for precise meaning). If we set $S_{\epsilon}^k(u)=\bigcap\limits_{r>0}S_{\epsilon,r}^k$ and $S^k(u)=\bigcup\limits_{\epsilon>0}S_{\epsilon}^k(u)$, then $S^k(u)$ coincides with the classical $k$-th stratum, which consists of points $x\in B_3(0)$ such that no tangent map of $u$ at $x$ is $(k+1)$-symmetric; see Lemma \ref{equi} below. The first main result of Breiner-Lamm \cite[Theorem 3.1]{Breiner-Lamm-2015} shows that if $u\colon B_3(x)\to N$ is a stationary biharmonic map in  $H_\Lambda^2(B_3(0),N)$ (see Definition \ref{dddd}, roughly speaking, with certain ``energy" no more than $\Lambda$), then for each $\epsilon>0$, there exists a constant $C=C(m,N,\Lambda,\epsilon)$ such that for each $0<r\leq 1$, we have
\begin{equation}\label{eq:B-L main result 1}
\text{Vol}\left(T_r(S_{\epsilon,r}^k(u))\cap B_1(x)\right)\leq Cr^{m-k-\epsilon},
\end{equation}
where $T_r(A)$ denotes the $r$-tubular neighborhood of a set $A$. This can be seen as a natural extension of the harmonic analogue obtained by Cheeger-Naber in \cite[Theorem 1.4]{Cheeger-Naber-2013-CPAM}.

On the other hand, in the case of harmonic maps, in another celebrated work, A. Naber and D. Valtorta \cite{Naber-V-2017} were able to improve the similar estimate as \eqref{eq:B-L main result 1} by removing the extra $\epsilon$ on the exponent of $r$ on the right hand side of \eqref{eq:B-L main result 1}. It is thus natural to ask
\medskip 

\textbf{Question. } \emph{Can we remove the extra $\epsilon$ on the exponent of $r$ on the right hand side of \eqref{eq:B-L main result 1}?  Namely, whether the following sharper estimate holds:}
\begin{align}\label{eq:expected}
	\mbox{\rm Vol}(T_r(S^k_{\ep,r}(u))\cap B_1(0))\leq C_\ep r^{m-k}.
\end{align}
\medskip 

The main motivation of this paper is to provide an affirmative answer to the above question. Our first main result reads as follows.
\begin{theorem}[Stratification of stationary biharmonic maps]\label{thm80}
Let $u\in H_\Lambda^2(B_3(0),N)$ be a stationary biharmonic map, where $H_\La$ is defined in Definition \ref{dddd}. Then for each $\ep>0$ there exists $C_\ep=C_\ep(m,N,\Lambda,\ep)$ such that for all $r\in(0,1]$,
\begin{align}\label{J}
\mbox{\rm Vol}(T_r(S^k_{\ep,r}(u))\cap B_1(0))\leq C_\ep r^{m-k}.
\end{align}
Consequently, for all $r\in(0,1]$,
 \begin{align}\label{Jk}
 \mbox{\rm Vol}(T_r(S^k_{\ep}(u))\cap B_1(0))\leq C_\ep r^{m-k}.
 \end{align}
 Moreover, for each $k$, $S^k_\ep(u)$ and $S^k(u)$ are $k$-rectifiable and  upper Ahlfors $k$-regular,  and  for $\HH^k$-almost every $x\in S^k(u)$, there exists a unique $k$-plane $V^k\su T_xM$ such that every tangent map of $u$ at $x$ is $k$-symmetric with respect to $V^k$.
\end{theorem}

Here a  subset $A\subset \R^n$ is said to be upper Ahlfors $k$-regular, if there is a constant $M>0$ such that 
\[\cal{H}^k(A\cap B_r(x))\le Mr^k\qquad \text{for all }\,x\in A \text{ and } 0<r<\text{diam}(A). \]
As a standard application of the improved volume estimate \eqref{J} in Theorem \ref{thm80}, we obtain optimal higher regularity for biharmonic maps. To state it, we need to recall the following definition of regularity scale.

\begin{definition}[Regularity scale]\label{pppv}
Let $u\colon B_2(0)\ri N$ be a map into manifold. For $x\in B_1(0)$, we denote the regularity scale $r_u(x)$ of $u$ at $x$ by
$$r_u(x):=\max\left\{0\leq r\leq 1:\sup_{B_r(x)}\left(r|\nabla u|+r^2|\nabla^2 u|+r^3|\nabla^3 u|+r^4|\nabla^4u|\right)\leq1\right\}.$$
\end{definition}
Note that if $r_u(x)=r>0$, then $u\in C^\infty(B_{{r}/{2}}(x))$. 

As an application of the volume estimate \eqref{eq:B-L main result 1}, C. Breiner and T. Lamm proved in their second main result \cite[Theorem 4.3 and Corollary 4.4]{Breiner-Lamm-2015} that if $u\colon B_4(0)\to N$ is minimizing biharmonic, then for all $\epsilon>0$, there exists $C=C(m,N,\Lambda,\epsilon)$ so that for each $0<r\leq 1$, 
\begin{equation}\label{eq:B-L main 2}
{\rm Vol}(\{x\in B_1(0):r_u(x)<r\})\leq Cr^{5-\epsilon}
\end{equation}
and so for all $1\leq p<5$ and $\ell=1,\cdots,4$, $|\nabla^{\ell}u|\in L^{p/\ell}(B_1(0))$. 

As an application of our improved volume estimate \eqref{J}, we are able to remove the extra $\epsilon$ on the exponent of $r$ in \eqref{eq:B-L main 2}. Consquently, we also obtain the optimal higher regularity for minimizing biharmonic maps: $|\nabla^{\ell}u|\in L^{5/\ell}_{weak}(B_1)$ for all $\ell=1,\cdots,4$. Here and after, $L^{p}_{weak}(B_1(0))$ denotes the space of weakly $L^p$-integrable functions on $B_1(0)$.

\begin{theorem}[Regularity estimates on minimizing biharmonic maps]\label{thm84} There exists a positive constant  $C=C(m,N,\La)$ such that, 
for any minimizing biharmonic map $u\in H_\Lambda^2(B_3(0),N)$, there holds
\begin{equation}\label{emhm}
\begin{aligned}
&\quad {\rm Vol}\Big(\{x\in B_1(0):r|\nabla u|+r^2|\nabla^2 u|+r^3|\nabla^3 u|+r^4|\nabla u|^4>1\}\Big)\\
&\leq{\rm Vol}(\{x\in B_1(0):r_u(x)<r\})\leq Cr^5.
\end{aligned}
\end{equation}
In particular, both $|\na^{\ell} u|$ and $r_u^{-\ell}$ have uniform bounds in $L^{5/\ell}_{weak}(B_1(0))$ for $\ell=1,2,3,4$. 

\end{theorem}

As was already observed in \cite[Remark 4.5]{Breiner-Lamm-2015}, the map $u\colon B_1^5\to S^4$ with $u(x)=\frac{x}{|x|}$ is minimizing biharmonic with $|\nabla u(x)|\approx 1/r$, which belongs to $L^5_{weak}$, but not in $L^p$ if $p>5$. This shows the sharpness of Theorem \ref{thm84} at least when $m=5$. 

As in the case of stationary/minimizing harmonic maps \cite{Naber-V-2017}, the preceeding regularity results can be further improved under extra topological assumptions. More precisely, we shall prove that a stationary/minimizing biharmonic maps have higher integrability if there are no smooth nonconstant ``biharmonic spheres" of certain dimensions, that is, biharmonic maps from $\R^{\ell+1}\backslash \{0\}$ to $N$. This observation  first appeared in the celebrated work of F.H. Lin \cite{Lin-1999-Annals} considering harmonic maps and then generalized to biharmonic maps by C. Breiner and T. Lamm in \cite[Theorem 4.3]{Breiner-Lamm-2015}, in which one of the crucial ingredients of the proof is the compactness theory established by  Scheven \cite{Scheven-2008-ACV} (see Theorems 1.5, 1.6 and Proposition 1.7 therein) and the useful auxiliary Lemma 4.7 of Breiner and Lamm \cite{Breiner-Lamm-2015}. We give a minor improvement of the corresponding result of Breiner-Lamm based on roughly the same argument.

\begin{theorem}[Improved estimates on biharmonic maps]\label{thm85}
Let $u\in H_\Lambda^2(B_3(0),N)$ be a stationary (minimizing resp.) biharmonic map. Assume that for some $k\geq4$ there exists no  nonconstant smooth 0-homogeneous stationary (minimizing resp.)  biharmonic maps $\R^{\ell+1}\backslash\{0\}\to N$ for all $4\leq \ell\leq k$. Then there exists a constant $C=C(m,K_N,\La)>0$ such that
\[\begin{aligned}
{\rm Vol}(\{x\in B_1(0)&:r|\nabla u|+r^2|\nabla^2 u|+r^3|\nabla^3 u|+r^4|\nabla^4 u|>1\})\\
&\leq{\rm Vol}(\{x\in B_1(0):r_u(x)<r\})\leq Cr^{2+k}.
\end{aligned}\]
In particular, both $|\na^\ell u|$ and $r_u^{-\ell}$ have uniform bounds in $L^{(2+k)/\ell}_{weak}(B_1(0))$ for $\ell=1,2,3,4$. 
\end{theorem}

A difficult but fundamental question would be \emph{whether one can obtain higher regularity estimates for stationary biharmonic maps $u\colon B^m\to N$, say $|\nabla^{\ell}u|\in L^{q/\ell}$ for some $4<q<5$ without any assumption on the target manifold}, as that of  minimizing biharmonic maps in Theorem \ref{thm84}. The question seems to be unknown even for stationary harmonic maps.

As in the case of harmonic maps, Theorem \ref{thm80} is the key result, and after it, all the later regularity theorems follow routinely. For the proof of Theorem \ref{thm80}, we shall follow closely the approach of A. Naber and D. Valtorta \cite{Naber-V-2017,Naber-V-2018}: establish the necessary quantitative $\epsilon$-regularity theorems and combine it with the Discrete Reifenberg Theorem proved in \cite{Naber-V-2017}, together with a delicate refined covering argument from \cite{Naber-V-2018}, to prove the improved volume estimate \eqref{J} and rectifiability of the quantitative stratum. Comparing with \cite{Naber-V-2017,Naber-V-2018}, the main difference lies in the monotonicity formula: the ``normalized energy" $\Phi_u$ contains not only the $L^2$-integral of $\Delta u$, but also some other term $H_u$ defined via integration on $\partial B_r$ (see Section \ref{sec:partial regularity} below for precise information). The earlier useful work of C. Scheven \cite{Scheven-2008-ACV} is not sufficient for application of the quantitative stratification method, and we need to develop further analysis of the defect measures for biharmonic maps to obtain symmetry of tangent maps and defect measures (see Proposition \ref{prop: tangent measur} below). This also makes the analysis on quantitative symmetries in Section \ref{sec:quantitative ep regularity} slightly more complicated than that developed for harmonic maps in \cite{Naber-V-2017,Naber-V-2018}. Once these were done, the covering argument from \cite{Naber-V-2018} can be relatively easily adapted to the biharmonic setting. For the convenience of readers, we made some effort to reformulate the covering lemmas (comparing with \cite{Naber-V-2018}) and slightly adjusted the proofs for the quantitative regularity (comparing with \cite{Breiner-Lamm-2015}) so that it is easier to understand the innovations introduced by Cheeger-Naber \cite{Cheeger-Naber-2013-CPAM} and Naber-Valtorta \cite{Naber-V-2017,Naber-V-2018}. With all these ingredients at hand, the proof of Theorem \ref{thm80} follows from a similar argument as that of Naber-Valtorta \cite{Naber-V-2017,Naber-V-2018} and Breiner-Lamm \cite{Breiner-Lamm-2015}.


\textbf{Notation.} The constant $C$ may be different from line to line in the proofs. A ball $B_r(0)\subset\mathbb{R}^m$ is  denoted by $B_r$ for simplicity.

\section{Partial regularity and defect measure}\label{sec:partial regularity}

In this section we collect some basic results for biharmonic maps.
\subsection{Monotonicity formula and partial regularity}

The following alternative describtion of weakly/stationary biharmonic maps are well-known; see e.g.~ \cite[Proposition 2.2]{Wang-2004-CPAM} and \cite[Lemma 1]{Angelsberg-2006}.
\begin{proposition}[\cite{Wang-2004-CPAM, Angelsberg-2006}]
	Suppose $u\in W^{2,2}(\Omega,N)$ is a weakly biharmonic map. Then there holds in the weak sense
	\begin{equation*}
	\De^2u=\De(A(u)(\nabla u,\nabla u))+2\nabla\cdot(\langle\nabla(P(u)),\De u\rangle)-\langle\De(P(u)),\De u\rangle.
	\end{equation*}
	Here 
	$P(u)=\nabla\pi_N(u)$, and $A(u)$ is the second fundamental form of $N$ at $u$.
	
	If, in addition, $u$ is stationary, then
	\begin{equation*}
	\int_\Omega4\De uD^2uD\psi+2\De uDu\De\psi-|\De u|^2\nabla\cdot\psi=0,
	\end{equation*}
	for any $\psi\in C_0^\infty(\Omega,\mathbb{R}^m)$, that is,
	\begin{equation}\label{1.5}
	\int_\Omega4u_{\al\al}\cdot u_{\beta\gamma}\psi_{\beta}^{\gamma}+2u_{\al\al}\cdot u_{\beta}\psi_{\gamma\gamma}^\beta-|\De u|^2\divg(\psi)=0.
	\end{equation}
	Here and hereafter we use the notation $u_\al=\frac{\pa u}{\pa x_\al}$ and Einstein's summation convention.
\end{proposition}
To continue, we introduce the following notations, whenever the integrals are well-defined:
\begin{eqnarray*}
& &H_u(a,r)\equiv r^{3-m}\int_{\partial B_r(a)}\Big(\partial_X|Du|^2+4|Du|^2-4\frac{|\partial_Xu|^2}{|x-a|^2}\Big)d\mathcal{H}^{m-1};\\
\label{eq: full normalized energy}
&&\Phi_u(a,r)\equiv r^{4-m}\int_{B_r(a)}|\De u|^2dx+H_u(a,r);\\
\label{eq: normalized energy}
&&G_u(a,r)\equiv r^{4-m}\int_{B_r(a)}|\nabla^2 u|^2+r^{-2}|\nabla u|^2dx;
\end{eqnarray*}
where $X=x-a$ and $\pa_Xu =(x-a) \cdot\na u$. 	Note that $\frac{|\partial_Xu|^2}{|x-a|^2}$ is the square norm of normal derivative of $u$ on the sphere $\pa B_r(a)$ and  $|Du|^2-\frac{|\partial_Xu|^2}{|x-a|^2}=|D_T u|^2$ is the square norm of tangent component of $Du$ on $\pa B_r(a)$.
\begin{remark} We will repeatedly use the translation and scaling invariance of biharmonic maps. More precisely, for any stationary biharmonic map $u$ on $\Om$, the function $u_{x,r}(y)=u(x+ry)$
	is a stationary biharmonic map on $r^{-1}(\Om-x)$. Moreover,
	\[
	H_{u_{x,r}}(0,1)=H_{u}(x,r),\quad	\Phi_{u_{x,r}}(0,1)=\Phi_{u}(x,r), \quad G_{u_{x,r}}(0,R)=G_u(x,rR).
	\]
\end{remark}

The following monotonicity formula plays a central role in the regularity issue of stationary biharmonic maps.
\begin{theorem}[Monotonicity formula,~\cite{Angelsberg-2006, Chang-W-Y-1999}]\label{thm1}
	Let $u\in W^{2,2}(B_R,N)$ be a stationary biharmonic map. Then, for any $a\in B_{R/2}$, $\Phi_u(a,r)$ and $H_u(a,r)$ are well defined for almost every $0<r<R/2$. Moreover,  for  almost every $ \ 0<\rho<r\leq {R}/{4}$, there holds
	\begin{eqnarray}\label{1.8}
	\Phi_u(a,r)-\Phi_u(a,\rho)=4\int_{B_r(a)\backslash B_\rho(a)}\left(\frac{|D\partial_Xu|^2}{|x-a|^{m-2}}+(m-2)\frac{|\partial_Xu|^2}{|x-a|^m}\right)dx.
	\end{eqnarray}	
	In particular, upon redefining $\Phi_u(a,\cdot)$ on a set of measure zero, we may assume \eqref{1.8} holds for all $0<\rho<r\leq R/4$ and thus $r\mapsto\Phi_u(a,r)$ is monotonically nondecreasing.
\end{theorem}
Using this monotonicity formula, we may define the density function $\Theta_u\colon \Omega\to [0,\infty)$ of a stationary biharmonic map $u\colon \Omega\to N$ as 
$$ \Theta_u(a):=\lim_{r\to 0} \Phi_u(a,r),\qquad a\in \Omega.$$ 
Letting $\rho\to 0$ in \eqref{1.8}, we obtain that for every $0<r<\dist(a,\pa\Omega)$, 
\[
\Phi_{u}(a,r)=\Theta_{u}(a)+4\int_{B_{r}(a)}\left(\frac{|D\pa_{X}u|^{2}}{|x-a|^{m-2}}+(m-2)\frac{|\pa_{X}u|^{2}}{|x-a|^{m}}\right)dx.
\]

The monotonicity formula implies that the gradient functions belong to certain Morrey space. More precisely, we have
\begin{lemma}\label{1.12} 	
	Suppose $u\in W^{2,2}(B_{4},N)$ is stationary biharmonic.	Then
	\begin{align*}
	\sup_{y\in B_1,\rho<1}\rho^{4-m}\int_{B_\rho(y)}\Big(|D^2u|^2+\rho^{-2}|D u|^2\Big)dx\leq C\int_{B_4}|D^2u|^2dx+\widetilde{C}
	\end{align*}
	for some constant $C$ depending only on $m$, and $\widetilde{C}$ depending only on $m$ and the diameter of $N\subset\mathbb{R}^L$.
\end{lemma}
\begin{proof}
	The estimate for $D^2 u$ can be found \cite[Lemma A.2]{Scheven-2008-ACV}. The estimate for $|D u|$ follows from an interpolation inequality of Gagliardo and Nirenberg in \cite{Nirenberg-1996}.
\end{proof}

Due to the monotonicity formula, the following partial regularity result holds, whose proof can be found in \cite[Theorem 1.3]{Wang-2004-CPAM}.
\begin{theorem}[\cite{Wang-2004-CPAM}]\label{thm3} 
	There exists a constant $\ep=\ep(m,N)>0$ such that for any stationary  biharmonic map $u\in W^{2,2}(\Omega,N)$, if  $G_u(y,R)\le \ep$,
	then $u\in C^\infty(B_{R/2}(y))$. Consequently, there exists a closed set  $\Sigma\subset\Omega$ given by
	\[\Sigma:=\left\{
	y \in \Omega : \liminf_{r \to 0} r^{4-m} \int_{B_r(y)}\left(\left|D^2 u\right|^2+r^{-2}|D u|^2\right) d x \geq \ep
	\right\},\]
	with  $\mathcal{H}^{m-4}(\Sigma)=0$, such that $u\in C^\infty(\Omega\backslash\Sigma,N)$.
\end{theorem}

Next, we define the singular set of a biharmonic map as follows.
\begin{definition}\label{singset}
	For any biharmonic map $u\in W^{2,2}(\Omega,N),$ we define its singular set $\sing(u)$ as
	$$\sing(u)=\{x\in\Omega:u\ \mbox{is not continuous at} \ x\},$$
	and its complementary set $\reg(u)$ as
	$$\reg(u):=\Omega\ba \sing(u).$$
\end{definition}

As in \cite{Breiner-Lamm-2015}, it is convenient to introduce the following function space.
\begin{definition}\label{dddd}
	For each $\Lambda>0$, the space $H_\Lambda^2(B_R(x),N)$ consists of all maps $u:B_R(x)\subset\mathbb{R}^m\rightarrow N$ such that
	\begin{equation*}
	G_u(x,R)\leq\Lambda
	\end{equation*} 
	and
	\begin{equation*}
	\Phi_u(x,R)\leq \Lambda.
	\end{equation*}
\end{definition}

For any map $u\in H^2(B_R(x),N)$, we can always choose a good radii  $R/4<r<R$ such that
	\begin{align*}
|H_u(x,r)|\leq CR^{4-m}\int_{B_{R}}(|\na^2 u|^2+R^{-2}|\na u|^2)dx 
\end{align*} 
for some $C>0$ depending only on $m$. This in turn implies that
	\begin{align*}
|\Phi_u(x,r)|\leq C_m R^{4-m}\int_{B_{R}}(|\na^2 u|^2+R^{-2}|\na u|^2)dx.
\end{align*} 
Hence for a  stationary biharmonic map $u\in H^2(B_R(x),N)$, if $G_u(x,R)\leq\Lambda$, then Lemma \ref{1.12} and the monotonicity formula implies that $G_u(x,R/4)\leq C\Lambda$ and $\Phi_u(x,R/4)\leq C\Lambda$. Therefore we have $u\in H^2_{C\Lambda}(B_{R/4},N)$.

\subsection{Defect measures}
In this subsection, we shall discuss the theory of defect measures for stationary biharmonic maps. In case of stationary harmonic maps, this was first considered in the seminal paper of F.H. Lin \cite{Lin-1999-Annals} and was extended to biharmonic maps by C. Scheven in \cite[Section 3]{Scheven-2008-ACV}. 
\begin{proposition}[Defect measures]\label{prop: defect measur} 
	Suppose that $\{u_{i}\}_{i\ge1}\subset H_{\Lambda}^{2}(B_{1},N)$
	is a sequence of stationary biharmonic maps
	satisfying
	\[
	u_{i}\wto u\quad\text{weakly in }H^{2}(B_{1}),\qquad u_{i}\to u\quad\text{strongly in }H^{1}(B_{1}) \text{ and almost everywhere }
	\]
	and
	\[
	|\De u_{i}|^{2}dx\wto\mu\qquad\text{in the sense of Radon measures.}
	\]
	Then the following conclusions hold:
	
	{\upshape (i)} There is a closed $(m-4)$-rectifiable set $\Sigma\subset B_{1}$ with
	${\cal H}^{m-4}(\Sigma)<\infty$, such that $u\in C_{\loc}^{\infty}(B_{1}\backslash\Sigma)$
	and
	\[
	u_{i}\to u\qquad\text{in } C_{\loc}^{2}(B_{1}\backslash\Sigma).
	\]
	
	{\upshape(ii)} $u$ is a weakly biharmonic map in $B_{1}$ with $\sing(u)\subset\Sigma$.
	
	{\upshape(iii)} There exist a nonnegative Radon measure $\nu$ with ${\rm spt}(\nu)\cup \sing(u)=\Sigma$ and a density function $\Theta_{\nu}\colon \Sigma\to [0,\infty)$  such that
	\[
	\mu=|\De u|^{2}dx+\nu,\qquad\text{and}\qquad \nu=\Theta_{\nu}{\cal H}^{m-4}|_{\Sigma}.
	\]
	Moreover, there are constants $C,c$ depending only on $m,n,N$ such that
	\[
	c\ep_{0}\le\Theta_{\nu}(x)\le C\Lambda\qquad \text{for } {\cal H}^{m-4}\text{-almost every } x\in \Sigma,
	\]
	where $\ep_0=\ep_0(m,N)$ is a constant from the $\ep$-regularity Theorem \ref{thm3}. 
	
	{\upshape(iv)} $u$ enjoys the unique continuation property in the sense that if there is another weakly harmonic map $v$ in $B_{1}$ such that $u=v$ almost everywhere on an open set, and $v$ is smooth away from a set $\Sigma'$ of finite $\mathcal{H}^{m-4}$-measure, then $u\equiv v$ on $B_{1}$. 
\end{proposition}

We call $\nu$ the \emph{defect measure} of the sequence $\{u_i\}$.

\begin{proof} 
	Following Scheven \cite[Equation (3.5) on page 63]{Scheven-2008-ACV}, we define the set $\Sigma$ by
	\[\Sigma:=\left\{a\in B_1: \liminf_{r\to 0} \left(r^{4-m}\int_{B_r(a)}|\De u|^2+r^2|\na u|^2dx+r^{4-m}\nu (B_r(a))\right)\ge \ep_0  \right\}.\]
	The first three assertions (i)--(iii) have been proved in \cite[Section 3]{Scheven-2008-ACV}, except the claim that  $\Sigma$ is $(m-4)$-rectifiable. This will be done with the aid of a deep theorem of Preiss \cite{Preiss-1987}. Indeed, by \cite[Theorem 3.4]{Scheven-2008-ACV}, we have
	\[c\ep_0\cal{H}^{m-4}|_\Sigma\le \nu|_{\overline{B_1}} \le C\Lambda\cal{H}^{m-4}|_{\Sigma}.\] 
	Thus the support spt($\nu$) of $\nu$ satisfies 
	\begin{equation}\label{eq:support of nu}
	\cal{H}^{m-4}(\Sigma\backslash{\rm spt}(\nu))=0.
	\end{equation}
	Thus,  for $\nu$-almost every $a\in B_1$, which is the same as $\cal{H}^{m-4}$-almost every $a\in B_1$, we have
	\[\limsup_{r\to 0} \frac{\nu (B_r(a))}{r^{m-4}} \le C\Lambda\limsup_{r\to 0} \frac{\cal{H}^{m-4}|_{\Sigma} (B_r(a))}{r^{m-4}} \le C\Lambda.\]
	On the other hand, by \cite[Equation (3.6)]{Scheven-2008-ACV}, there holds
	\[ \liminf_{r\to 0} \frac{\nu (B_r(a))}{r^{m-4}} \ge \ep_0, \qquad \text{for $\nu$-almost every } a\in B_1. \]
	Combining these two estimates, we obtain  
	 \[0<\limsup_{r\to 0} \frac{\nu (B_r(a))}{r^{m-4}} \le C \liminf_{r\to 0} \frac{\nu (B_r(a))}{r^{m-4}}<\infty \]
	for $\nu$-almost every  $a\in B_1$.  Therefore, by \cite[Theorem 1.2]{DeLellis-book}, $\nu$ is ($m-4$)-rectifiable, and consequently by \eqref{eq:support of nu}, $\Sigma$ is also ($m-4$)-rectifiable. 	
	
	 Thus we only need to prove assertion (iv). 	
	Note that if $u$ and $v$ are both smooth, then the claim follows from \cite[Theorem 1.5]{Branding-Oniciuc-2019}. In the general case, let $\Sigma$ be the closed set as in (i). Then $u$ and $v$ are smooth outside the closed set $\Sigma_0:=\Sigma\cup \Sigma'$. This closed set is at most ($m-4$)-dimensional and thus non-disconnecting. Applying \cite[Theorem 1.5]{Branding-Oniciuc-2019} again gives $u=v$ almost everywhere on $B_1\backslash \Sigma_0^c$ and thus on the whole domain. 
\end{proof}




\begin{remark}  By the unique continuation property above, if $u$ is a weakly biharmonic map as above, and is radially invariant on $B_r(x)\backslash B_\rho(x)$ for some $x\in \mathbb{R}^m$ and some $0<\rho<r<1$, then we can extend $u$ to the whole space $\mathbb{R}^m$ as a radially invariant (with respect to $x$) weakly biharmonic map.
\end{remark}

Now we use defect measure theory to explore properties of tangent maps and tangent measures of stationary biharmonic maps. A basic observation is  that if $u\in H^2(\Omega, N)$ is a stationary biharmonic map, then so is the rescaled map $u_{a,r}:= u(a+r\cdot)$  for any $a\in \Omega$ and $r>0$.
\begin{definition}
	We say that $v\in H_{\loc}^{2}(\mathbb{R}^m,N)$ is a tangent map of a stationary harmonic map $u\in H^2(\Omega, N)$ at the point $a\in \Omega$, if there is a sequence $r_j\to 0$ such that $u_{a,r_j}\equiv u(a+r_j\cdot)\wto v$ in $H_{\loc}^{2}(\mathbb{R}^m,N)$.
\end{definition}

Based on the theory of defect measures (see Proposition \ref{prop: defect measur}), we are able to deduce the following basic results for tangent maps of biharmonic maps.
\begin{proposition} \label{prop: tangent measur}
	Let $u\in H^2(B_1, N)$ be a stationary biharmonic map. Suppose there is a sequence $r_j\to 0$ such that
	\[
	u_{j}\wto v\quad\text{in }H^{2}_\loc(\mathbb{R}^m, N),\qquad u_{j}\to v\quad\text{in }H_\loc^{1}(\mathbb{R}^m, N) \text{ and almost everywhere},
	\] and
	\[	|\De u_{j}|^{2}dx\wto\mu\qquad\text{in the sense of Radon measures,}	\]
	where $u_{j}=u_{a,r_j}:=u(a+r_j\cdot)$ for $a\in B_1$. 	Then
	
	{\upshape(i)} There is a closed $(m-4)$-rectifiable set $\Sigma\subset \mathbb{R}^m$ with locally finite 
	${\cal H}^{m-4}$-measure, such that $v\in C_{\loc}^{\infty}(\mathbb{R}^m\backslash \Sigma)$
	and
	\[
	u_{j}\to v\qquad\text{in } C_{\loc}^{2}(\mathbb{R}^m\backslash\Sigma).
	\]
	
	{\upshape(ii) (Symmetry of tangent map)} $v$ is a  weakly biharmonic map  with ${\rm sing}(v)\subset\Sigma$, and moreover, it is 0-homogeneous with respect to the origin, i.e.
	$$v(x)=v(\lambda x),\qquad\text{for all } x\in\R^m,\lambda>0.$$
	
	{\upshape(iii)} There exist a nonnegative Radon measure $\nu$ with ${\rm spt}(\nu)\cup \sing(v)=\Sigma$ and a density function $\Theta_{\nu}\colon  \R^m \to [0,\infty)$  such that
	\[
	\mu=|\De v|^{2}dx+\nu\qquad\text{and}\qquad \nu=\Theta_{\nu}{\cal H}^{m-4}|_{\Sigma}
	\]
	Moreover, there are nonnegative constants $C,c$ depending only on $m,n,N$ such that
	\[
	c\ep_{0}\le\Theta_{\nu}(x)\le C\Lambda\qquad{\cal H}^{m-4}\text{-almost everywhere on }\Sigma,
	\]
	where $\ep_0$ is given in {\upshape (iii)} of Proposition \ref{prop: defect measur}.
	
	{\upshape(iv) (Symmetry of tangent measure)} $\nu$ is a cone measure in the sense that for any $\la>0$ and measurable set $A\subset\mathbb{R}^m$, there holds
	\[\la^{4-m}\nu(\la A)=\nu(A).\]
	Consequently, due to the homogeneity of $v$ in assertion {\upshape(ii)}, the measure $\mu$ is a cone measure as well. 
\end{proposition}
We call $\mu$ or $(v,\nu)$ as in Proposition \ref{prop: tangent measur} a \emph{tangent measure} of $u$ at the point $a\in B_1$.

\begin{proof}
	
	Assertions (i)--(iii) are direct consequences of Proposition  \ref{prop: defect measur}, except the conclusion concerning the homogeneity of $v$. This fact is surely well-known for experts. However, as it is the fundation of our later analysis, we sketch the proof here for convenience.
	
	Using the monotonicity formula, we find that
	\begin{equation}\label{sinaf}
	\Phi_{u_{j}}(0,r)-\Phi_{u_{j}}(0,\rho)\geq4(m-2)\int_{B_r(0)\backslash B_\rho(0)}\frac{|\pa_X u_{j}|^2}{|x|^m}dx,
	\end{equation}
	where $\pa_Xu:=x\cdot\na u$.
	On the other hand, for any $r>0$ we have $\Phi_{u_{j}}(0,r)=\Phi_u(a,r_jr)$ for all $j\ge 1$.
	Hence by the definition of the density function $\Theta_u$ it follows that 
	$$\lim_{j\to\oo}\Phi_{u_{j}}(0,r)=\lim_{j\to\oo}\Phi_{u }(a,r_jr)=\lim_{r\to 0}\Phi_{u }(a,r)=\Theta_{u}(a).$$
	Therefore, by letting $j\to\infty$ in \eqref{sinaf} and using the strong convergence $u_{j}\to  v$ in $W^{1,2}_{\loc}$ (up to a subsequence) we get
	\begin{equation*}
	\int_{B_r(0)\backslash B_\rho(0)} \frac{|\pa_Xv|^2}{|x|^m} dx=0.
	\end{equation*}
	This implies that $v$ is radially invariant in $B_r(0)\backslash B_\rho(0)$. Since $r,\rho$ are arbitrary, this implies that $v$ is radially invariant with respect to the origin.

	Next we prove the symmetry of $\nu$ claimed in assertion (iv) via the method of Lin \cite{Lin-1999-Annals} for stationary harmonic maps.
	Given any  $\psi\in C^2_c(\mathbb{R}^m)$, we aim  to show that
	\begin{equation}\label{eq: scaling invariance of tangent measure}
	\int_{\mathbb{R}^m}\Big((m-4)\psi_\la+x\cdot\na\psi_{\la}\Big) d\mu=0,
	\end{equation} where $\psi_\lambda=\psi({\cdot}/{\la})$ for  $\lambda>0$.
	Equivalently, this means that
	\[	
	\frac{d}{d\la}	\int_{\mathbb{R}^m} \psi d\mu_{0,\la}=\frac{d}{d\la}	\int_{\mathbb{R}^m} \la^{4-m} \psi_\la d\mu\equiv 0 \qquad\text{for all }\la>0\,
	\] 
	where $\mu_{0,\la}$ is defined by $\mu_{0,\la}(A)=\la^{4-m}\mu(\la A)$. Once  the above claim were true, then $\mu_{0,\la}\equiv\mu_{0,1}= \mu$ for all $\la>0$, showing the cone measure property of $\nu$. Consequently, $\nu=\mu-|\De v|^2dx$ is a cone measure as well, as being the difference of two cone measures.
	
	To prove \eqref{eq: scaling invariance of tangent measure}, substituting  $\zeta(x)=\psi_\la(x)x$ into the stationary equation \eqref{1.5} yields
	\begin{equation}\label{q1.5a}
	0=\int_{\mathbb{R}^m}(4\De u_i(u_i)_{\al\beta}\pa_\al(\psi_\la(x) x^\beta)+2\De u_i(u_i)_\beta\De(\psi_\la(x) x^\beta)
	-|\De u_i|^2\divg(\psi_\la(x) x)).
	\end{equation}
	By the weak convergence $|\Delta u_i|^2dx\wto \mu$, we have
	\begin{align}\label{q1.5da}
	-\lim_{i\to\oo}\int_{\mathbb{R}^m} |\De u_i|^2\divg(\psi_\la(x) x)=-\int_{\mathbb{R}^m}(\na\psi_\la(x)\cdot x+m\psi_\la)d\mu.
	\end{align}
	Since $u_i\to v$ strongly in $H^1_{\loc}$ and $\De u_i\wto \De v$ weakly in $L^2_\loc$, we infer that
	\begin{equation}\label{eq: 2rd term}
	\begin{aligned}
	\lim_{i\to \infty} \int_{\mathbb{R}^m}2\De u_i(u_i)_\beta\De(\psi_\la(x) x^\beta)dx&= \int_{\mathbb{R}^m}2\De v\cdot v_\beta\De(\psi_\la(x) x^\beta)dx\\
	&=\int_{\R^m}2\Delta v\Delta\psi_{\lambda}\nabla v\cdot x+4\Delta v\nabla v\cdot \nabla \psi_{\lambda} dx.
	\end{aligned}
	\end{equation}
	The first term in \eqref{q1.5a} can be computed directly to give
	\begin{align}\label{q1.0a}
	\int_{\mathbb{R}^m}4\De u_i(u_i)_{\al\beta}\pa_\al(\psi_\la(x) x^\beta)
	=\int_{\mathbb{R}^m}4\De u_i(u_i)_{\al\beta}\pa_\al\psi_\la(x) x^\beta +\int_{\mathbb{R}^m}4|\De u_i|^2\psi_\la.
	\end{align}
	The weak convergence $|\Delta u_i|^2dx\wto \mu$ implies that
	\begin{equation}\label{eq: item 1-2}\lim_{i\to \infty}\int_{\mathbb{R}^m}4|\De u_i|^2\psi_\la dx=\int_{\mathbb{R}^m}4\psi_\la d\mu.\end{equation}
	We are left to dispose the first term of the right hand side of \eqref{q1.0a}. This will be done by applying the monotonicity formula.
	
	Recall that for any $0<\rho<r$, there holds
	\begin{align}\label{q61}
	\Phi_{u_i}(0,r)-\Phi_{u_i}(0,\rho)=4\int_{B_r\backslash B_{\rho}}\left((m-2)\frac{|x\cdot\na u_i|^2}{|x|^m}+\frac{|D\left(x\cdot\na u_i\right)|^2}{|x|^{m-2}}\right)dx,
	\end{align}
	and
	$$\lim_{i\ri\oo}\Phi_{u_i}(0,r)=\lim_{i\ri\oo}\Phi_{u }(a,rr_i)=\Theta_{u}(a).$$
	Since both terms on the right hand side of the above equality are nonnegative,  by letting $\rho \to 0$ first in \eqref{q61}  and then letting $i\to \infty$ we deduce that
	\begin{align}\label{q6cb}
	\lim_{i\ri\oo}\int_{B_r}{\sum_{\beta} |\pa_\beta \left(x\cdot\na u_i\right)|^2}dx=\lim_{i\ri\oo}\int_{B_r}{|D\left(x\cdot\na u_i\right)|^2}dx=0.
	\end{align}   
	On the other hand, since $\pa_\beta \left(x\cdot\na u_i\right)= (u_i)_{\beta}+x^\al\cdot (u_i)_{\al\beta}$ and $(u_i)_{\beta}\to v_{\beta}$ strongly in $L^2_\loc$, we conclude from \eqref{q6cb} that for all $1\le \beta\le m$, 
	\begin{align}\label{q6cb-2}
	x^\al\cdot (u_i)_{\al\beta} \to -v_{\beta} \qquad \text{in } L^2_{\loc}(\R^m).
	\end{align}
	Hence
	\begin{equation}\label{eq: item 1-1}
	\int_{\mathbb{R}^m}4\De u_i(u_i)_{\al\beta}\pa_\al\psi_\la(x) x^\beta dx\to -\int_{\mathbb{R}^m}4\De v \nabla v\cdot \nabla \psi_\la dx.
	\end{equation}
	Combining \eqref{eq: item 1-1}, \eqref{eq: item 1-2} and  \eqref{q1.0a} gives
	\begin{equation}\label{eq: item 1}
	\int_{\mathbb{R}^m}4\De u_i(u_i)_{\al\beta}\pa_\al(\psi_\la(x) x^\beta)dx\to \int_{\mathbb{R}^m}4\psi_\la d\mu-4\De v \nabla v\cdot \nabla \psi_\la dx.
	\end{equation}
	Hence, combining \eqref{q1.5a}, \eqref{q1.5da}, \eqref{eq: 2rd term} and \eqref{eq: item 1} we obtain
	\[	\int_{\mathbb{R}^m}\Big((m-4)\psi_\la+x\cdot\na\psi_{\la}\Big) d\mu+\int_{\mathbb{R}^m} 2\Delta v(x\cdot \na v)\De \psi_\lambda=0.\]
	As $v$ is radially invariant, we have $x\cdot \na v=0$.	The desired equation~\eqref{eq: scaling invariance of tangent measure} follows from \eqref{eq: item 1}, \eqref{eq: 2rd term} and \eqref{q1.5da}. The proof of assertion (iv) and thus also Proposition \ref{prop: tangent measur} is complete.
\end{proof}


\section{Quantitative $\ep$-regularity theorem}\label{sec:quantitative ep regularity}

In the previous section, we have seen that tangent maps/measures of stationary harmonic maps enjoy certain symmetry and one can follow the classical approach to develop a stratification theory of the singular set for biharmonic maps based on the symmetry of tangent maps/measures (see Appendix \ref{sec:appendix}). On the other hand, for the proof of Theorem \ref{thm80}, we shall adapt the quantitative stratification method for harmonic maps developed by Naber-Valtorta \cite{Naber-V-2017}. Towards this, we introduce quantitative symmetry and study properties of quantitative stratum. The proofs of these results are in spirit similar to that of Naber-Valtorta \cite{Naber-V-2017,Naber-V-2018} but some care need to be paid due to the complication of monotonicity formula. For the convenience of readers, we will give full details for the proofs as much as possible. 

\subsection{Quantitative symmetry and cone splitting principle}
\begin{definition}[Symmetry]\label{def:symmetry}
Given a measurable map $h\colon \mathbb{R}^m\to N$ we say that

{\upshape(1)} $h$ is 0-homogeneous or 0-symmetric with respect to  point $p$ if $h(p+\lambda v)=h(p+v)$ for all $\lambda>0$ and $v\in\mathbb{R}^m$.

{\upshape(2)} $h$ is $k$-symmetric if it is 0-homogeneous with respect to the origin,  and is translation invariant with respect to a $k$-dimensional subspace $V\subset\mathbb{R}^m$, i.e.,
$$h(x+v)=h(x) \qquad \text{for all } x\in\mathbb{R}^m,\, v\in V.$$
\end{definition}
If $u\in C^1(\mathbb{R}^m,N)$, then $u$ is  0-homogeneous at $x=p$ if and only if $\partial_{r_p}h=0$, where $r_p$ is the radial direction centered at $p$; and $u$ is translation invariant with respect to a $k$-dimensional subspace $V\subset\mathbb{R}^m$ if and only if $\pa_v u=0$ for all $v\in V$.

\begin{definition}[Quantitative  symmetry]\label{def:quantitative symmetry}
Given a map $u\in L^2(\Omega,N)$, $\ep>0$ and nonnegative integer $k$,  we say that  $u$ is $(k,\ep)$-symmetric  on $B_r(x)\subset \Omega$, or simply $B_r(x)\subset\Omega$ is $(k,\ep)$-symmetric,  if there exists some $k$-symmetric function $h\colon \R^m\to N$ such that
$$\medint_{B_r(x)}|u(y)-h(y-x)|^2\leq\ep.$$
\end{definition}
Note that  $u$ is $(k,\ep)$-symmetric  on $B_r(x)$ if and only if the scaled map  $u_{x,r}(y)=u(x+ry)$ is $(k,\ep)$-symmetric  on $B_1(0)$.

Given the definition of quantitative symmetry, we can  introduce a quantitative stratification for points of a function  according to how much it is symmetric around those points.
\begin{definition}[Quantitative stratification]\label{ddd}
For any map $u\in L^2(\Omega,N)$, $r,\ep>0$ and $k\in\{0,1,\cdots,m\}$, we define the $k$-th quantitative singular stratum $S_{\ep,r}^k(u)\subset \Omega$ as 
$$S^k_{\ep,r}(u)\equiv\big\{x\in\Omega: u \text{ is not } (k+1,\ep)\text{-symmetric on }B_s(x) \text{ for any } r\le s<1\big\}.$$
Furthermore, we set
\[
S^k_\ep(u):=\bigcap_{r>0}S^k_{\ep,r}(u)
\quad\text{ and }\quad 
S^k(u)=\bigcup_{\ep>0}S^k_\epsilon(u).
\]
\end{definition}

Note that by definition,
$$k'\leq k\text{ or }\ep'\geq\ep\text{ or } r'\leq r\ \Longrightarrow \ S^{k'}_{\ep',r'}(u)\subseteq S^k_{\ep,r}(u).$$
In particular, we have
\[S^0(u)\subset S^1(u)\subset\cdots\subset S^m(u)=\Omega.\]
 
The following lemma shows that $S^k(u)$ is well adapted to stationary biharmonic maps.
 
\begin{lemma}\label{equi} Suppose $u\in H^2(\Om,N)$ is a stationary biharmonic map. Then
	$$S^k(u)=\{x\in \Omega: \text{no tangent maps of } u \text{ at } x \mbox{ is } (k+1)\mbox{-symmetric}\}.$$
	Consequently, we have
	\[S^0(u)\subset S^1(u)\subset\cdots\subset S^{m-1}(u)\subset\sing(u).\]
\end{lemma}
\begin{proof} 
	For the moment we write 
	\[\Sigma^k(u)=\{x\in \Omega: \text{no tangent maps of } u \text{ at } x \mbox{ is } (k+1)\mbox{-symmetric}\}.\]
	
	Suppose $x\in S^k(u)$. Then $x\in S^k_\ep(u)$ for some $\ep>0$. Thus, for any $(k+1)$-symmetric map $h\in L^2$ and $r>0$, we have
	$$\medint_{B_1(0)}|u_{x,r}-h|^2dx\geq\ep.$$
	If $v$ is a tangent map of $u$ at $x$, there exists a sequence $r_i\to 0$ such that $u_{x,r_i}\rightharpoonup v$ in $W^{2,2}_{\loc}(\R^m,N)$. Then it follows that
	$$\medint_{B_1(0)}|v-h|^2dx\geq\ep,$$
	which implies that $v$ is not $(k+1)$-symmetric.  Hence
	$S^k(u)\subset \Sigma^k(u)$.
	
	For the reverse direction, suppose $x\notin S^k(u)$. Then there exist sequences of positive numbers $r_i>0$ and $(k+1)$-symmetric maps $h_i\colon \R^m\to N$ such that
	$$\medint_{B_1(0)}|u_{x,r_i}-h_i|^2dx\leq i^{-1}.$$
	A further extraction of subsequences (still denoted by $i $) leads to $u_{x,r_i}\rightharpoonup v$ in $W^{2,2}$ and $h_i\rightharpoonup h$ in $L^2(B_1(0))$. Then  by the weak lower semi-continuity of $L^2$-norm, we obtain
	\begin{align*}
	\int_{B_1(0)}|v-h|^2dx\leq\liminf_{i\ri\oo}\int_{B_1(0)}|u_{y,r_i}-h_i|^2dx=0.
	\end{align*}
	{Moreover, by the compactness of symmetric maps, $h$ is $(k+1)$-symmetric.}

	If $r_i\ri0$,  then $v$ is a tangent map and thus is $(k+1)$-symmetric, which shows that $x\not\in  \Sigma^k(u)$. If $\lim_{i\ri\oo}r_i=r>0$, then there exists a positive radius $\delta\leq r$ such that $u(z)=h(z-x)$ for almost every $z\in B_\delta(x)$. Indeed, notice that
	\[
	\medint_{B_{r_i}(x)}|u(z)-h_i(z-x)|^2dz=\medint_{B_1(0)}|u_{x,r_i}(z)-h_i(z)|^2dz\leq i^{-1}
	\]
	and so we conclude
	\[
	\medint_{B_\delta(x)}|u(z)-h(z-x)|^2dz=0,
	\]
	by sending $i\to \infty$ and using the lower semi-continuity of $L^2$-norm. On the other hand, the previous indentity implies that all tangent maps of $u$ at $x$ are $(k+1)$-symmetric and thus $x\notin \Sigma^k(u)$.  	
\end{proof}

Note that it might be possible that there is a constant tangent map for a stationary biharmonic map $u$ at a singular point $x$. 
\begin{remark}
	Due to the strong convergence in \cite[Theorem 1.5]{Scheven-2008-ACV}, if $u$ is a minimizing biharmonic map, then we  have the equality $S^{m-5}(u)=\sing(u)$.
\end{remark}

Next, we establish  quantitative  symmetry and cone splitting principle  for  stationary biharmonic maps $u$ via $\Phi_u$.

Note that if  $u$ is a stationary biharmonic map and $\Phi_u(a,1)=\Phi_u(a,1/2)$, then by the monotonicity formula \eqref{1.8} and the unique continuation property (see Proposition \ref{prop: defect measur} (iv)), $u$ must be $0$-symmetric with respect to $x=a$. The following proposition gives a  quantitative version of this property.
\begin{proposition}\label{pro1}	For any $\ep>0,$ there exists $\delta_0=\delta_0(m,N,\Lambda,\ep)$ such that, if  $u\in H^2_\La(B_{8},N)$ is a stationary biharmonic map with
	$$\Phi_u(x,r)-\Phi_u(x,r/2)<\delta_0$$
	for some $x\in B_1$ and $0<r<1$, then $u$ is $(0,\ep)$-symmetric on $B_r(x)$.
\end{proposition}
\begin{proof}
	Suppose, on  the contrary,  that for some $\ep>0$ there is a sequence of  stationary biharmonic  maps $u_i\in H^2_{\Lambda}(B_8)$ which satisfy
	$$\Phi_{u_i}(x_i,r_i)-\Phi_{u_i}(x_i,r_i/2)<i^{-1},$$  but is not $(0,\ep)$-symmetric on $B_{r_i}(x_i)$.
	Let $\bar{u}_i(y)=u_i(x_i+r_iy)$. Then  
	$$\Phi_{\bar{u}_i}(0,s)=\Phi_{u_i}(x_i,r_is).$$ 
	By Lemma \ref{1.12}, $\{\bar{u}_i\}_{i\geq 1}$ are uniformly bounded in $W^{2,2}(B_1).$ So up to  a subsequence, we may  assume that $\bar{u}_i\rightharpoonup v$ weakly in $W^{2,2}(B_1)$ and $\bar{u}_i\to v$ strongly in $W^{1,2}$. Then, by Proposition \ref{prop: defect measur},   $v$ is weakly biharmonic in $B_1$ with $\dim_{\mathcal{H}}\sing(v)\le m-4$.
	
	Now using the monotonicity formula \eqref{1.8}, for all $i\in \N$, we have
	\[
	\int_{B_1\backslash B_{1/2}}{|y\cdot\nabla \bar{u}_i|^2}dy \le C(\Phi_{\bar{u}_i}(0,1)-\Phi_{\bar{u}_i}(0,1/2))=C(\Phi_{u_i}(x_i,r_i)-\Phi_{u_i}(x_i,r_i/2)).
	\]
    Thus, by sending $i\to \infty$ and using the strong convergence of $\bar{u}_i\rightarrow v$ in $W^{1,2}(B_1)$, we deduce that $v$ is radially invariant on $B_1(0)\backslash B_{1/2}(0)$,  and thus is homogeneous on $\R^m$ by the unique continuation property (see Proposition \ref{prop: defect measur}). Finally, since $\bar{u}_i\rightarrow v$ strongly in $L^2$, it follows that
	\[ 
	\medint_{B_{r_i}(x_i)} |u_i-v((\cdot-x_i)/r_i)|^2=\medint_{B_{1}}|\bar{u}_i-v|^2\to 0\qquad \text{as }i\to \infty, 
	\] 
	contradicting with the assumption that $u_i $ is not  $(0,\ep)$-symmetric on $B_{r_i}(x_i)$. The proof is thus complete. 
\end{proof}

Next, we extend Proposition \ref{pro1} to higher order symmetric cases. For this, we first introduce the notion of  quantitative frame for affine subspaces as in \cite{Naber-V-2018}.
\begin{definition}
	Let $\{y_i\}_{i=0}^k \subset B_1(0)$ and $\rho>0$. We say that these points $\rho$-effectively span a $k$-dimensional affine subspace if for all $i=1,\cdots,k$,
	$$\dist(y_i,y_0+\Span\{y_1-y_0,\cdots,y_{i-1}-y_0\})\geq2\rho.$$
	More generally, a set $F\subset B_1(0)$ is said to $\rho$-effectively span a $k$-dimensional affine subspace, if there exist points $\{y_i\}_{i=0}^k\subset F$ which $\rho$-effectively spans a $k$-dimensional affine subspace.
\end{definition}

\begin{remark}\label{rem4.5}
	Such a set $\{y_i\}_{i=0}^k$ is also called $\rho$-independent, and produces a {quantitative frame}   $\{y_i-y_0\}_{i=1}^k$. The advantage of quantitative frame is twofold (see the comments right below \cite[Definition 28]{Naber-V-2018}):
	
	 {\upshape(i)} If $\{y_i\}_{i=0}^k$ $\rho$-effectively span a $k$-dimensional affine subspace, then for every point 	$x\in y_0+{\rm span}\{y_1-y_0,\cdots,y_k-y_0\}$,
	there exists a unique set of numbers $\{\alpha_i\}_{i=1}^k$ such that
	$$x=y_0+\sum_{i=1}^k\alpha_i(y_i-y_0) \quad \text{with}\quad  |\alpha_i|\leq C(m,\rho)|x-y_0|.$$
	
	{\upshape(ii)} Quantitative frame is stable under limiting process: if  $\{y_i^j\}_{i=0}^k$ $\rho$-effectively spans a $k$-dimensional affine subspace for all $j\ge 1$, and $y_i^j\to y_i$ as $j\to \infty$, then $\{y_i\}_{i=0}^k$ also $\rho$-effectively spans a $k$-dimensional affine subspace.
\end{remark}

Before stating the next result, we point out an useful observation, called the \emph{cone splitting principle} \cite{Cheeger-Naber-2013-CPAM}: If $u\colon \mathbb{R}^n\rightarrow\mathbb{R}^m$ is $k$-symmetric with respect to a $k$-plane $V^k,$ and is $0$-symmetric at a point $z\not\in V^k$, then $u$ is $(k+1)$-symmetric with respect to  the $(k+1)$-plane $V^{k+1}:=\Span\{z,V^k\}.$ 

The following result, which can be viewed as a quantitative version of the above cone splitting principle, extends Proposition \ref{pro1} to higher order quantitative symmetry.
\begin{proposition}[Quantitative cone splitting]\label{5.2}
For any $\ep,\rho>0$, there exists $\delta_1=\delta_1(m,\Lambda,N,\ep,\rho)>0$ such that, if $u\in H^2_{\Lambda}(B_8,N)$ is a stationary biharmonic map, and  for some $\{x_j\}_{j=0}^k\subset B_1(0)$ there holds
		
	{\upshape(a)} $\Phi_u(x_j,r)-\Phi_u(x_j,r/2)<\delta_1$ for some $0<r<1$ and for all $j=0,\cdots,k$,
	
	{\upshape(b)} $\{x_j\} \ \rho$-effectively spans a $k$-dimensional affine plane $V$,
	
	\noindent then $u$ is $(k,\ep)$-symmetric on  $B_1(0)$.
\end{proposition}
\begin{proof}
	Suppose by contradiction that there are a sequence of  stationary biharmonic maps $u^i\colon B_8(0)\rightarrow N$ with
	$\|u^i\|_{W^{2,2}(B_8,N)}\leq\Lambda$,  $\ep_0,\rho_0>0$, $\{x_0^i,\cdots,x_k^i\}$ and a sequence of $k$-dimensional affine plane $V^i$ such that $\{x_0^i,\cdots,x_k^i\}$  $\rho_0$-effectively spans $V^i$ and $\Phi_{u^i}(x^i_j,r)-\Phi_{u^i}(x^i_j,r/2)<\delta_i$, where $\delta_i\to 0$ is chosen according to Proposition \ref{pro1}, so that $u^i$ is $(0,i^{-1})$-symmetric at $B_{i^{-1}}(x_j^i)$ for each $j$. But $u^i$ is not $(k,\ep_0)$-symmetric in $B_1(0)$ for all $i\geq1.$  
	
	Assume $u^i\rightharpoonup v$ weakly in $W^{2,2}$ for some biharmonic map $v\in W^{2,2}(B_2(0))$ and $u^i\rightarrow v$ strongly in $W^{1,2}(B_2(0)).$ Assume also $x_j^i\rightarrow x_j\in B_1(0)$ and $V^i\ri V$ as $i\rightarrow\infty.$
	Applying the monotonicity formula to each $u^i$, we obtain
	$$\int_{B_r(x^i_j)\backslash B_{r/2}(x_j^i)}|(x-x^i_j)\cdot\na u^i|^2\leq C\left(\Phi_{u^i}(x^i_j,r)-\Phi_{u^i}(x^i_j,r/2)\right)<C\delta_i.$$
	Sending $i\ri\oo$, we get
	$$\int_{B_r(x_j)\backslash B_{r/2}(x_j)}|(x-x_j)\cdot\na v|^2=0.$$
	This shows that $v$ is radially symmetric on ${B_r(x_j)\backslash B_{r/2}(x_j)}$ for each $j\in \{0,\cdots,k\}$. The unique continuation property for biharmonic map (see \cite{Branding-Oniciuc-2019}) and the cone splitting principle imply $v$ is symmetric with respect to the $k$-dimensional affine plane spanned by $\{x_j\}_{j=0}^k$. However, this contradicts with the assumption that $u^i$ is not $(k,\ep_0)$-symmetric in $B_1(0)$ and the strong convergence of $u^i\to v$ in $W^{1,2}(B_1).$
\end{proof}

\subsection{Properties of quantitative stratum} The next proposition shows that under certain assumptions the singular set satisfies a  one-side Reifenberg approximating property.
\begin{proposition}\label{pro32} 
	For any $\ep,\rho>0$, there exists $\delta_2=\delta_2(m,N,\Lambda,\ep,\rho)>0$ such that the following holds: If $u\in H^2_{\Lambda}(B_8,N)$ is a stationary biharmonic map, and the set $$F:=\{y\in B_2: \Phi_u(y,2)-\Phi_u(y,\rho)<\delta_2\}$$ $\rho$-effectively spans a $k$-dimensional affine plane $V$, then $$S^k_{\ep,\delta_2}(u)\cap B_1\subset B_{2\rho}(V).$$
\end{proposition}

The proof of Proposition \ref{pro32} relies on the following lemma, which says that if the map $u$ inside $B_1$ is almost translation invariant along some $(k+1)$-dimensional subspace, then  $u$ is almost $(k+1)$-symmetric in a neighborhood of every point in $B_{1/2}$.
\begin{lemma}\label{5.3}
	For any $\ep>0$, there exists $\delta_3=\delta_3(m,N,\Lambda,\ep)>0$ such that if  $u\in H^2_{\Lambda}(B_8,N)$ is a stationary biharmonic map satisfying
	\begin{align}\label{81}
	\int_{B_1(0)}|P\cdot\nabla u|^2<\delta_3
	\end{align}
	for some $(k+1)$-dimensional subspace $P$, then 
	$S^k_{\ep,\overline{r}}(u)\cap B_{1/2}=\emptyset$ for  $\overline{r}=\delta_3^{\frac{1}{2(m-2)}}$. In particular, $0\not\in S^k_{\ep,\overline{r}}(u)$.
\end{lemma}

Here and after, $|P\cdot\nabla u|^2=\sum\limits_{i=1}^{k+1}|\nabla_{e_i}u|^2$ with $\{e_i\}_{i=1}^{k+1}$ being an  orthonormal basis of $P$.

\begin{proof}	
	First we claim that there is a constant $C_2(m,n,\Lambda)>0$ such that for every $x\in B_{1/2}(0)$, there exists $r_x\in[\overline{r},1/2]$  such that
	\begin{align}\label{80}
	\Phi_u(x,r_x)-\Phi_u(x,r_x/2)<\frac{C_2(m,n,\Lambda)}{-\log\delta_3}.
	\end{align}
	Indeed, if this is not true, then for  $u\in H^2_\Lambda(B_8)$ we may assume by choosing a good radii that $\Phi_u(x,1/2)\le C_1(m,n,\Lambda)$. Then 
	\begin{align*}
	C_1(m,n,\Lambda)\geq\Phi_u(x,1/2)&\geq\sum_{i=1}^{-\log\overline{r}+1}(\Phi_u(x,2^{-i})-\Phi_u(x,2^{-i-1}))\\
	&\geq c(m)C_2(m,n,\Lambda),
	\end{align*}
	which is impossible if we take $$C_2(m,n,\Lambda)=\frac{2C_1(m,n,\Lambda)}{{c}(m)\Lambda}.$$ This proves the claim.

	We now use a contradiction argument to prove our claim. Suppose there exist some $\ep>0$ and a contradicting sequence of stationary biharmonic maps  $u_i\in H^2_\Lambda(B_8)$, together with subspaces $P_i$, $\de_{3,i}\to 0$, $x_i\in B_{1/2}(0)$ and $r_i\in[\bar{r}_i,1]$, such that $u_i$ is not $(k+1,\ep)$-symmetric on $B_{r_i}(x_i)$, but \eqref{80} holds for $x=x_i$, where $\bar{r}_i=\de_{3,i}^{1/2(m-2)}$ and $r_i=r_{x_i}$. Note also that by the definition of $r_x$, we have
	\begin{align}\label{79}
	r_i^{2-m}\int_{B_{r_i}(x_i)}|P_i\cdot\nabla u_i|^2<r_i^{2-m}\delta_{3,i}\leq\de_{3,i}^{1/2}.
	\end{align} 
		Using a simple rotation, we may assume that the $(k+1)$-dimensional subspaces $P_i$ are fixed by $P$, i.e., $P_i=P$ for all $i$. 
	
	Let $v_i(x)=u_i(x_i+r_ix)$. Then we may assume that $v_i$ converges weakly in $H^2$ and strongly in $H^1$ to some weakly biharmonic map $v$. It follows that $v$ is 0-symmetric by unique continuation and is invariant with respect to the $(k+1)$-dimensional subspace $P$ by \eqref{79}: 
	\[
	0=\lim_{i\to \infty}r_i^{2-m}\int_{B_{r_i}(x_i)}|P\cdot\nabla u_i|^2=\lim_{i\to \infty}\int_{B_1}|P\cdot \nabla v_i|^2=\int_{B_1}|P\cdot \nabla v|^2.
	\]
	Since $v_i\to v$ in $L^2$, this implies that $v_i$ is $(k+1,\ep)$-symmetric on $B_{1}$ if $i\gg 1$, or equivalently, $u_i$ is $(k+1,\ep)$-symmetric with $B_{r_i}(x_i)$. This is clearly a contradiction.
\end{proof}

Now we can prove Proposition \ref{pro32}.
\begin{proof}[Proof of Proposition \ref{pro32}]
Suppose $\{y_{j}\}_{j=0}^{k}\subset F$ is a $\rho$-independent frame
that spans the affine subspace $V$, i.e., 
\[
V=y_{0}+{\rm span}\{y_{i}-y_{0}\}_{i=1}^{k}.
\]
Let $x_{0}\in B_{1}\backslash B_{2\rho}(V)$ and $\de_{2}>0$ to be
determined later. We need to prove that $x_{0}\not\in S_{\ep,\de_{2}}^{k}(u)\cap B_{1}$. 

By the definition of $F$, there holds $\Phi_{u}(y_{i},2)-\Phi_{u}(y_{i},\rho)<\de_{2}$
for all $0\le i\le k$. Note that, for $0<r\ll\rho$, we have $B_{r}(x_{0})\subset B_{2}(y_{i})\backslash B_{\rho}(y_{i})$
for every $i\le k$. Thus, for each $0\le i\le k$,
\[
\int_{B_{r}(x_{0})}|(z-y_{i})\cdot\na u|^{2}dz\le\int_{B_{2}(y_{i})\backslash B_{\rho}(y_{i})}|(z-y_{i})\cdot\na u|^{2}dz\le C\de_{2}.
\]
by the monotonicity formula \eqref{1.8} of $\Phi_{u}$. Consequently,
\[
\begin{aligned}\int_{B_{r}(x_{0})}|(y_{i}-y_{0})\cdot\na u|^{2} & \le2\int_{B_{2}(y_{i})\backslash B_{\rho}(y_{i})}|(z-y_{i})\cdot\na u|^{2}+2\int_{B_{2}(y_{i})\backslash B_{\rho}(y_{i})}|(z-y_{0})\cdot\na u|^{2}\\
& \le2C\de_{2}.
\end{aligned}
\]
Since $\{y_{j}\}_{j=0}^{k}$ is $\rho$-independent, we conclude that
\begin{equation}
\int_{B_{r}(x_{0})}|\hat{V}\cdot\na u|^{2}\le C\de_{2}\label{eq: k-small}
\end{equation}
for some constant $C=C(m,\rho)>0$, where $\hat{V}={\rm span}\{y_{i}-y_{0}\}_{i=1}^{k}$.

On the other hand, for each $z\in B_{r}(x_{0})\subset B_{2}\backslash B_{\rho}(V)$,
let 
\[
\pi_{V}(z)=y_{0}+\sum_{i=1}^{k}\al_{i}(z)(y_{i}-y_{0})
\]
be the orthogonal projection of $z$ in $V$. Then $|z-\pi_{V}(z)|\ge\rho$,
$|\al_{i}(z)|\le C(m,\rho)$ and
\[
\int_{B_{r}(x_{0})}|(z-\pi_{V}(z))\cdot\na u|^{2}dz\le C\sum_{i=0}^{k}\int_{B_{2}(y_{i})\backslash B_{\rho}(y_{i})}|(z-y_{i})\cdot\na u|^{2}\le C(m,\rho)\de_{2}.
\]
Thus, by setting $h(z)=(z-\pi_{V}(z))/|z-\pi_{V}(z)|$, it follows
that
\[
\begin{aligned}\int_{B_{r}(x_{0})}|h(x_{0})\cdot\na u|^{2} & \le\int_{B_{r}(x_{0})}|h(z)\cdot\na u|^{2}+\int_{B_{r}(x_{0})}|(h(z)-h(x_{0}))\cdot\na u|^{2}\\
& \le C\de_{2}+Cr^{2}\int_{B_{r}(x_{0})}|\na u|^{2}\\
& \le C\de_{2}+C(m,\rho,\La)r^{m}.
\end{aligned}
\]
Now we choose $r=r(m,\La,\rho)\ll\rho$ such that $C(m,\rho,\La)r^{m}\le C\de_{2}$.
It gives
\[
\int_{B_{r}(x_{0})}|h(x_{0})\cdot\na u|^{2}\le2C\de_{2}.
\]
Together with \eqref{eq: k-small}, for the $(k+1)$-dimensional subspace $P:=\hat{V}\oplus \R h(x_0)$, we find that
\[
\int_{B_{r}(x_{0})}|P\cdot\na u|^{2}\le C(m,\rho,\La)\de_{2}.
\]
Now choosing $\de_{2}=\de_{2}(m,\rho,\La)$ sufficiently small so
that we can apply Lemma \ref{5.3} to conclude that $u$ is $(k+1,\ep)$-symmetric
in $B_{r}(x_{0})$. This completes the proof of Proposition \ref{pro32}.
\end{proof} 

The next result shows that $\Phi_u$ remains almost constant on all pinched points.
\begin{lemma}\label{5.5}
 Let $\rho,\eta>0$ be fixed and  $u\in H^2_{\Lambda}(B_8,N)$  a stationary biharmonic map. Let  
 $$E=\sup_{y\in B_2(0)}\Phi_u(y,1).$$ 
 Then there exists $\de_4=\de_4(m,N,\Lambda,\rho,\eta)>0$ such that, if the set 
 $$F:=\{y\in B_2:\Phi_u(y,\rho)>E-\delta_4\}$$
  $\rho$-effectively spans a $k$-dimensional affine subspace $L\subset\mathbb{R}^m$, then
 $$\Phi_u(x,\rho)\ge E-\eta \qquad \text{for all }  x\in L\cap B_2(0).$$
\end{lemma}
\begin{proof}
We prove this lemma by contradiction. Suppose $\{u_i\}\subset H^2_\La(B_8,N)$ is a sequence of stationary biharmonic maps with
$$\sup_{y\in B_2(0)}\Phi_{u_i}(y,1)\leq E.$$
For each $i\ge 1$, we may assume that $F_i:=\{y\in B_2:\Phi_{u_i}(y,\rho)>E-i^{-1}\}$ contains a subset $\{y^i_j\}_{j=0}^k$ spanning $\rho$-effectively  a $k$-dimemsional affine subspace $L_i\subset\mathbb{R}^m$, and there exists $x_i\in L_i\cap B_2(0)$ such that
\begin{align}\label{fanz}
\Phi_{u_i}(x_i,\rho)\leq E-\eta.
\end{align}
It follows from the assumption that
\begin{equation}\label{eq:lemma 312}
\Phi_{u_i}(y^i_j,1)-\Phi_{u_i}(y^i_j,\rho)<1/i,\qquad \text{for all }\,i\ge 1,\, 0\le j\le k. 
\end{equation}
Without loss of generality, we assume that for each $0\le j\le k$
\[y^i_j\to y_j\quad\text{and}\quad x_i\to x \qquad \text{as }i\to \infty\]
and that the sequence  $L_i$ of $k$-dimensional affine subspaces converges to a $k$-dimensional affine subspace $L$ passing through  $x$.

By Proposition \ref{prop: defect measur}, there exists a Radon measure $\mu=|\De v|^2dx+\nu$, where $v$ is a weakly biharmonic map and $\nu$ is the defect measure such that up to a subsequence $|\De u_i|^2dx\wto \mu$. Using the monotonicity formula \eqref{1.8} and the stationarity of $u_i$, adapting the argument in the proof of Proposition \ref{prop: tangent measur},  we find that for any fixed $y$, the function
$$r\mapsto \Phi_{\mu}(y,r):=\Phi_v(y,r)+r^{4-m}\nu(B_r(y))$$ is monotonically nondecreasing (upon modifying $\Phi_v(y,\cdot)$ on a set of measure zero). Thus sending $i\to \infty$ in \eqref{eq:lemma 312} and \eqref{fanz}, we obtain
\[\Phi_{\mu}(y_j,1)=\Phi_{\mu}(y_j,\rho)=E,\qquad \text{for all }\, 0\le j\le k\]
and \begin{equation}\label{eq: loss of energy}
\Phi_\mu(x,\rho)\le E-\eta.
\end{equation}
Moreover, as that of  Proposition \ref{prop: tangent measur}, we know that  $\mu,v,\nu$ are  translation invariant along $L$. Hence $\Phi_\mu|_{L\cap B_2(0)}\equiv \Phi_{\mu}(y_j,1)=E$, which clearly contradicts with \eqref{eq: loss of energy}. The proof is complete.
\end{proof}

The following technical result shows that the almost symmetry is preserved under certain pinching condition.
\begin{lemma}\label{5.6}
 For any  $\ep,\rho>0$, there exists $\de_5=\de_5(m,N,\Lambda,\rho,\ep)>0$ satisfying the following property. Suppose $u\in H^2_{\Lambda}(B_8,N)$ is a stationary biharmonic map and
 $\Phi_u(0,1)-\Phi_u(0,1/2)<\de_5$. If there is a point $y\in B_3$ such that

 {\upshape 1)}. $\Phi_u(y,1)-\Phi_u(y,1/2)<\de_5$,

 {\upshape 2)}.  $u$ is not $(k+1,\ep)$-symmetric on $B_r(y)$ for some $r\in[\rho,2]$,\\
then $u$  is not $(k+1,\ep/2)$-symmetric on $B_r(0)$. 
\end{lemma}
In particular, under condition {\upshape 1)} of Lemma \ref{5.6},  $y\in S^k_{\ep,\rho}(u)\cap B_3\Rightarrow 0\in S^k_{\ep/2,\rho}(u)$.

\begin{proof}
	Suppose by contradiction that there is a sequence $\{u_i\}_{i\in \N}$ of stationary biharmonic maps  satisfying $\Phi_{u_i}(0,1)-\Phi_{u_i}(0,1/2)\leq i^{-1},$ and there exists a sequence $y_i\in B_3$ such that $\Phi_{u_i}(y_i,1)-\Phi_{u_i}(y_i,1/2)\leq i^{-1}$ and that for each $i\in \N$, $u_i$ is not $(k+1,\ep)$-symmetric on $B_r(y_i)$, but it is $(k+1,\ep/2)$-symmetric on $B_r(0)$. 
	
	Then, there exists a sequence of $(k+1)$-symmetric maps $h_i\colon B_{4}(0)\ri N$ such that
	$$\medint_{B_r(0)}|u_i-h_i|^2\leq\ep/2,\quad i=1,2\cdots.$$
	Up to a subsequence if necessary, we may assume that $y_i\ri y\in \overline{B}_3(0),u_i\rightharpoonup v$ in $W^{2,2},$  and $u_i\ri v$ in $W^{1,2}$ for some $v\in W^{2,2}$, and $h_i\rii h$ in $L^2$. By the unique continuation for biharmonic maps, we see that $u$ is a weakly biharmonic map and is homogeneous with respect to the origin and $y.$ By the compactness of symmetric functions, $h$ is $(k+1)$-symmetric. By the property of weak convergence we have
	$$\medint_{B_r(0)}|v-h|^2\leq\limsup_{i\ri\oo}\medint_{B_r(0)}|v-h_i|^2\leq2\ep/3.$$
	If $y=0,$ then since $N$ is compact, $\|u_i\|_{L^\oo}<+\infty$, and we get
	\begin{align*}
		\lim_{i\ri\oo}\medint_{B_r(y_i)}&|u_i(x)-h(x-y_i)|^2\leq\lim_{i\ri\oo}\medint_{B_r(y_i)}|u_i(x)-h(x)|^2dx\\
		&+\lim_{i\ri\oo}\medint_{B_r(y_i)}|h(x)-h(x-y_i)|^2dx\leq2\ep/3.
	\end{align*}
	If $y\neq0,$ then by the invariance of $v$ with respect to $\R y$, we get
	\begin{align*}
		\lim_{i\ri\oo}&\medint_{B_r(y_i)}|u_i(x)-h(x-y_i)|^2=\lim_{i\ri\oo}\medint_{B_r(0)}|u_i(y_i+x)-h(x)|^2dx\\
		&\leq 2\lim_{i\ri\oo}\medint_{B_r(0)}|u_i(y_i+x)-v(y+x)|^2dx+2\medint_{B_r(0)}|v(y+x)-h(x)|^2dx\\
		&\leq 4\ep/3.
	\end{align*}
	In both cases, we arrive at a contradiction.
\end{proof}

\begin{remark}\label{rmk:on Lemma 313}
In Section \ref{sec:covering lemma}, we shall repeatedly use the following variant of Lemma \ref{5.6}: Suppose $u\in H^2_{\Lambda}(B_3(x),N)$ is a stationary biharmonic map and $\Phi_u(x,1)-\Phi_u(x,1/2)<\de_5$. If there is some $y\in B_3(x)$ with $\Phi_u(y,1)-\Phi_u(y,1/2)<\de_5$, then $y\in S^k_{\ep,\rho}(u)\cap B_3(x)\Rightarrow x\in S^k_{\ep/2,\rho}(u)$.
\end{remark}

\section{Reifenberg theorems and estimates of Jones' number}
In this section, we recall the Reifenberg type results obtained by Naber-Valtorta~\cite{Naber-V-2017}. 
Following the presentation in \cite{Naber-V-2018}, we first recall Jones'  number  $\beta_2$ which quantifies how close the support of a measure $\mu$ is to a $k$-dimensional affine subspace.

\begin{definition}\label{aaa}
Let $\mu$ be a nonnegative Radon measure on $B_3$. Fix  $k\in\N$. The $k$-dimensional Jones' $\beta_2$ number is defined as
\begin{align*}
\beta^k_{2,\mu}(x,r)^2=\inf_{V\su\R^m}\int_{B_r(x)}\frac{d^2(y,V)}{r^2}\frac{d\mu(y)}{r^k},\qquad \text{when }\,B_r(x)\subset B_3,
\end{align*}
where the infimum is taken over  all  $k$-dimensional affine subspaces $V$.
\end{definition}

We now state two versions of the quantitative Reifenberg theorems from \cite{Naber-V-2017}.
 \begin{theorem}[{Discrete-Reifenberg, \cite[Theorem 3.4]{Naber-V-2017}}]\label{rre}
 There exist two constants $\de_6=\de_6(m)$ and $C_R(m)$ such that the following property holds. Let $\{B_{r_x}(x)\}_{x\in \mathcal{C}}\su B_2(0)\su\R^m$ be a family of pairwise disjoint balls with centers in $\mathcal{C}\subset B_1(0)$ and let $\mu\equiv\sum_{x\in \mathcal{C}}\omega_kr_x^k\de_x$ be the associated measure. If for every ball $B_r(x)\su B_2,$ there holds
 \begin{align}\label{eq: multiscale appro-1}
\int_{B_r(x)}\left(\int_0^r\beta^k_{2,\mu}(y,s)^2\frac{ds}{s}\right)d\mu(y)<\de_6^2r^k,
\end{align}
then we have the uniform estimate
 \begin{align*}
\sum_{x\in \mathcal{C}}r^k_x< C_R(m).
\end{align*}
 \end{theorem}

\begin{theorem} [{Rectifiable-Reifenberg, \cite[Theorem 3.3]{Naber-V-2017}}]\label{rree}
 There exist constants $\de_7=\de_7(m)$ and $C=C(m)$ satisfying the following property. Assume that $S\subset B_2\su\R^m$ is $\HH^k$-measurable, and for each $B_r(x)\subset B_2$ there holds
 \begin{align}\label{eq: multiscale appro-2}
\int_{S\cap B_r(x)}\left(\int_0^r\beta^k_{2,\HH^k|_S}(y,s)^2\frac{ds}{s}\right)d\HH^k(y)<\de_7^2r^k.
\end{align}
Then $S\cap B_1$ is  $k$-rectifiable, and   $\HH^k(S\cap B_r(x))\leq C r^k$ for each $x\in S\cap B_1$.
 \end{theorem}

The two conditions \eqref{eq: multiscale appro-1} \eqref{eq: multiscale appro-2} are usually called multiscale approximation conditions. How to control the number $\beta_2$ is the key problem in application. 
Similar to the case of stationary harmonic maps \cite[Theorem 7.1]{Naber-V-2017},  we need to establish an $L^2$ subspace approximation theorem.

For $x\in B_1(0)$ and $r>0,$ we denote for brevity
 \begin{align*}
 W_r(x):= W_{r,10r}(x)=\int_{B_{10r}(x)\ba B_r(x)}\frac{|(y-x)\cdot\na u(y)|^2}{|y-x|^m}dy\geq0.
\end{align*}

\begin{theorem}\label{ji}
Fix $\ep>0,0<r\leq1$ and $x\in B_1(0).$ Then there exist $C(m,N,\Lambda,\ep)>0$ and $\de_8>0,$ such that, if $u\in H^2_\La(B_{10},N)$ is a stationary biharmonic map, and is $(0,\de_8)$-symmetric on $B_{10r}(x)$ but not $(k+1,\ep)$-symmetric, then for any nonnegative finite measure $\mu$ on $B_r(0)$ we have
\begin{align}\label{po}
\beta^k_{2,\mu}(x,r)^2\leq Cr^{-k}\int_{B_r(x)}W_r(y)d\mu(y).
\end{align}
\end{theorem}
In order to prove this theorem, we can assume that $\mu$ is  a probability measure supported on  $B_1(0)$. Let $x_{cm}$ be the mass center of $\mu$ in  $B_1(0)$, i.e.,  
$
x_{cm}\equiv\int xd\mu(x).
$
Here and later the integral domain is  $B_1(0)$ and is omitted for brevity. The second moment $Q=Q(\mu)$ of $\mu$ is the symmetric bilinear form  defined by
\begin{align*}
	Q(v,w):=\int[(x-x_{cm})\cdot v][(x-x_{cm})\cdot w]d\mu(x),\qquad \text{for all }\ v,w\in \R^m.
\end{align*}
Let $\lambda_1(\mu)\ge \cdots\ge \lam_m(\mu)$ be nonincreasing eigenvalues of $Q(\mu)$ and $v_1(\mu),\cdots,v_m(\mu)$ be the associated eigenvectors. Then we have
\begin{align}\label{il}
	Q(v_k)=\lam_kv_k=\int[(x-x_{cm})\cdot v_k][(x-x_{cm}) ]d\mu(x).
\end{align}
We can also describe the eigenvalues by variational method, that is
\begin{align*}
	\lam_1=\lam_1(\mu):=\max_{|v|^2=1}\int|(x-x_{cm})\cdot v|^2d\mu(x).
\end{align*}
Let $v_1=v_1(\mu)$ is any unit vector achieving such maximum. By induction we have
\begin{align*}
	\lam_{k+1}=\lam_{k+1}(\mu):=\max\left\{\int|(x-x_{cm})\cdot v|^2d\mu(x):{|v|^2=1} ,\ v\cdot v_i=0,\ \text{for all }i\leq k\right\},
\end{align*}
 and $v_{k+1}=v_{k+1}(\mu)$ is any unit vector achieving such maximum. Note that by definition of  $v_k$, it is not difficult to show that  $V_k=x_{cm}+\mbox{span}\{v_1,\cdots,v_k\}$ is the $k$-dimensional affine subspace achieving the minimum in the definition of $\beta_2$; see \cite[Remark 49]{Naber-V-2018}. Moreover,
\begin{align*}
\beta_{2,\mu}^k(0,1)^2=\int d^2(x,V_k)d\mu(x)=\lam_{k+1}(\mu)+\cdots+\lam_m(\mu).
\end{align*}
In fact, the second equality follows from $d^2(x,V_k)=\sum_{i=k+1}^{m}((x-x_{cm})\cdot v_i)^2$
and the definition of $\la_k$.
The following property gives the relationship of $\lam_k,v_k$ and $W_1$.
\begin{proposition}[{\cite[Proposition 50]{Naber-V-2018}}]\label{hj}
	Let  $\mu$ be a probability measure in $B_1(0)$ and $u\in H^1(B_{10}(0), N)$. Let $\lam_k,\ v_k$ be defined  as  above. Then there exists $C(m)>0$ such that
	\begin{align*}
		\lam_k\int_{A_{3,4}}|v_k\cdot\na u(z)|^2dz\leq C(m)\int_{B_1} W_1(x)d\mu(x), \qquad \text{for all }\,k\ge 1,
	\end{align*}
	where $A_{3,4}=B_4(0)\backslash B_3(0)$.
\end{proposition}


To further estimate the left hand side of the estimate in the above proposition, we need the following reversed result of  Lemma \ref{5.3}.
\begin{lemma}\label{lem: non higher order sym}
	For any $\ep>0$, there exists  $\de_9=\de_9(m,N,\La,\ep)>0$ satisfying the following property. For any stationary biharmonic map $u\in H^2_\La(B_8,N)$, if it is $(0,\de_9)$-symmetric on $B_1(0)$ but not $(k+1,\ep)$-symmetric, then 	
	\begin{align}\label{81}
	\int_{A_{3,4}}|P\cdot\nabla u|^2>\delta_9
	\end{align}
	for every $(k+1)$-dimensional subspace $P$.
\end{lemma}
\begin{proof}
	Assume by contradiction that there is a sequence of stationary biharmonic maps $u_i\in H^2_\La(B_8,N)$ such that $u_i$ is $(0,i^{-1})$-symmetric  on $B_1(0)$ but not $(k+1,\ep_0)$-symmetric for some $\ep_0>0$. Moreover, after an orthogonal transformation, there is a  $(k+1)$-dimensional subspace $P$ such that
	\begin{align}\label{999}
	\int_{A_{3,4}}|P\cdot\nabla u_i|^2\leq i^{-1},\ i=1,2,\cdots.
	\end{align}
	By extracting a subsequence, we can assume $u_i\rightharpoonup v$ in $W^{2,2}(B_8)$ and $u_i\ri v$ in $W^{1,2}(B_8)$ for some 0-homogeneous biharmonic map $v\in W^{2,2}(B_8)$ and 
	\begin{align*}
	\int_{A_{3,4}}|P\cdot\nabla v|^2=0.
	\end{align*}
	By the unique continuation property we know 
	$	\int_{B_4}|P\cdot\nabla v|^2=0.$ 
	Thus $v$  is $(k+1)$-symmetric on $B_1(0)$. Consequently $u_i$ is $(k+1,\ep_0)$-symmetric for $i\gg 1$, since $u_i$ converges strongly to $v$ in $L^2$. This is a contradiction and thus the proof is complete. 
\end{proof}
 Now we can prove Theorem \ref{ji}.
\begin{proof}[Proof of Theorem \ref{ji}]  By scaling, we may assume $\mu(B_1(0))=1.$ Since $\lambda_k$ is nonincreasing,
\begin{align}\label{pl}
	\beta_{2,\mu}^k(0,1)^2=\lam_{k+1}+\cdots+\lam_m\leq(m-k)\lam_{k+1}.
\end{align}
Thus it suffices to estimate $\la_{k+1}$. 
For each $j=1,\cdots,k+1$, by Proposition \ref{hj}, we have
\begin{align*}
	\sum_{j=1}^{k+1}\lam_j\int_{A_{3,4}(0)}|\na u(z)\cdot v_j|^2dz\leq(k+1)C\int W_1d\mu(x).
\end{align*}
Let $V^{k+1}=\mbox{span}(v_1,\cdots,v_{k+1}),$ then
\begin{align*}
	\lam_{k+1}\int_{A_{3,4}(0)}|V^{k+1}\cdot\na u(z)|^2dz=\lam_{k+1}\sum_{j=1}^{k+1}\int_{A_{3,4}(0)}|\na u(z)\cdot v_j|^2dz\leq C\int W_1d\mu(x).
\end{align*}
Let $\de_9$ be the number defined as in Lemma \ref{lem: non higher order sym} and set $\delta_8=\delta_9$. By assumption $B_{10}(0)$ is $(0,\de_8)$-symmetric but not $(k+1,\ep)$-symmetric.  Then by Lemma \ref{lem: non higher order sym}, there holds 
\begin{align*}
	\int_{A_{3,4}(0)}|V^{k+1}\cdot\na u(z)|^2dz\geq \de_8,
\end{align*}
and so
\begin{align*}
	\de_8\lam_{k+1}\leq\lam_{k+1}\int_{A_{3,4}(0)}|V^{k+1}\cdot\na u(z)|^2dz \leq C\int W_1d\mu(x).
\end{align*}
Since $\de_8=\de_9$ depends only on $m,N,\La,\ep$, by \eqref{pl}, we conclude 
\begin{align*}
	\beta_{2,\mu}^k(0,1)^2\leq C(m,N,\La,\ep)\int W_1d\mu(x).
\end{align*}
This completes the proof.
\end{proof}


\section{Covering lemma}\label{sec:covering lemma}
In this section we establish a covering lemma which aims to decrease the energy of biharmonic maps at singular points. We shall follow the presentation of Naber-Valtorta \cite[Section 6.2]{Naber-V-2018}, but with some further simplifications. Our main result of this section is to prove the following covering result. 

\begin{lemma}[Main covering Lemma]\label{lemma:main covering lemma}
	Let $u\in H_{\La}^{2}(B_3,N)$ be a stationary biharmonic map.
	Fix any $\ep>0$ and $0<r<R\leq1$. Then there exist constants $\de=\de(m,N,\La,\ep)>0$ and $C(m)$ with the following property. For any
	${S}\su S^k_{\ep,\de r},$ there exists a finite covering of $S\cap B_R(0)$ such that
	\begin{align}\label{pla}
		S\cap B_R(0)\su \bigcup_{x\in\3}B_{r_x}(x) \quad\text{with}\quad  r_x\leq r \quad\text{and}\quad  \sum_{x\in\3}r_x^k\leq C(m)R^k.
	\end{align}
	Moreover, the balls in $\{B_{r_x/5}(x)\}_{x\in \3}$ are pairwise disjoint and $\cal{C}\subset 	S\cap B_R(0) $. 
\end{lemma}

To prove Lemma \ref{lemma:main covering lemma}, we divide the proof into two sublemmata. 

\subsection{The first covering lemma}
We first establish the following lemma.
\begin{lemma}[Covering Lemma I]\label{oooa}
	Let $u\in H_{\La}^{2}(B_3,N)$ be a stationary biharmonic map.
	Fix any $\ep>0$, $0<\rho=\rho(m)\leq100^{-1}$ and $0<r<R\leq1$. Then there exist constants $\de=\de(m,N,\La,\rho,\ep)>0$ and $C_1(m)$ with the following property: For any
	${S}\su S^k_{\ep,\de r}$ with
	$$E:=\sup\limits_{x\in B_{2R}(0)\cap S}\Phi_u(x,R)\leq\Lambda, $$ 
	there exists a finite covering of $S\cap B_R(0)$ such that
	\begin{align}\label{pla}
		S\cap B_R(0)\su \bigcup_{x\in\3}B_{r_x}(x), \quad\text{with}\quad  r_x\leq r \quad\text{and}\quad  \sum_{x\in\3}r_x^k\leq C_1(m)R^k.
	\end{align}
	Moreover, for each $x\in\3,$ one of the following conditions  is satisfied:
	\begin{itemize}
		\item[i)] $r_x=r$;
		\item[ii)] the set of points $F_x:=\{y\in S\cap B_{2r_x}(x):\Phi_u(y,\rho r_x/10)>E-\de\}$ is contained in $B_{\rho r_x/5}(L_x)\cap B_{2r_x}(x),$ where $L_x$ is some  $(k-1)$-dimensional affine subspace.
	\end{itemize}
\end{lemma}
\begin{proof}
	By the scaling invariance property of $\Phi$, we may assume $R=1$. For cleaness of the presentation, we also assume in the following proof that
	\begin{align}\label{332}
		r=\rho^{{l}}, \quad \rho=2^{-a},\quad a,\ {l}\in\N.
	\end{align}
	
	We divide the proof into two main  steps.
	
	\textbf{Step 1. (Inductive covering)} Let $\eta>0$ be a constant to be determined later. Define the initial set
	\[
	F^0=\{y\in B_2(0)\cap{S}:\Phi_u(y,\rho/10)\geq E-\de\}.
	\]
	 If there is a $(k-1)$-dimensional affine subspace $L_0$ such that $F^0\su B_{\rho/5}(L_0),$ then we say that $B_1(0)$ is a \emph{good ball}.  If there is a  $k$-dimensional affine subspace $V_0$,  which is ${\rho}/{10}$-effectively spanned by some $\{y_j\}_{j=0}^k\su F^0$, then we call it a \emph{bad ball}.
	
	If $B_1(0)$ is a good ball, then our claim clearly holds. In the latter case, if $\de\leq\de_2$ sufficiently small, then we can apply Proposition \ref{pro32} to $B_1(0)$ to obtain
	\begin{equation}
	S^k_{\ep,\de r}\cap B_1(0)\su S^k_{\ep,\de}\cap B_1(0)\su B_{\rho/5}(V_0).
	\end{equation}
	Choose a finite cover of $B_{\rho/5}(V_0)\cap B_1$ by balls $\{B_\rho(x)\}_{x\in\3^1}$ with $\3^1\su V_0\cap B_1(0)$ such that if $x\neq y$ and $x,y\in\3^1$, then $B_{\rho/5}(x)\cap B_{\rho/5}(y)=\emptyset.$
	Note that by Lemma \ref{5.5} we have for all $x\in\3^1\su V_0\cap B_1(0),$
	$$\Phi_u(x,\rho/10)\geq E-\eta,$$
	as long as $\de$ is sufficiently small. Under the same smallness assumption, Lemma \ref{5.6} (or Remark \ref{rmk:on Lemma 313}) implies that for each $x\in \3^1$ we have $x\in S^k_{\ep/2,\rho}.$	
This gives the covering 
	$$S\cap B_1(0)\su\bigcup_{x\in \3^1}B_{\rho}(x)=\bigcup_{x\in \3^1_b}B_{\rho}(x)\cup\bigcup_{\3^1_g}B_{\rho}(x)\equiv B_{\rho}(\3^1_b)\bigcup B_{\rho}(\3^1_g),$$
	where
	\begin{align*}
		\3^1_g:=\left\{x\in\3^1:\ F_x^1 \su B_{ \rho^2/5}(L_x^1) \
		\mbox{\rm for some}\ (k-1)\mbox{\rm -dimensional affine subspace}\ L_x^1\right\}
	\end{align*}
	with
	\begin{align*}
		F_x^1=\left\{y\in S\cap B_{2\rho}(x):\Phi_u(y, \rho^2/10)\geq E-\de\right\}
	\end{align*}
	and
	\begin{align*}
		\3^1_b:=\{x\in&\3^1:\ F_x^1\ (\rho^2/10)\mbox{-}\mbox{\rm effectively spans a}\ k\mbox{\rm -dimensional affine subspace}\ L_x^{1'}\}.
	\end{align*}
	
	Consider any bad ball $B_\rho(x)$.
    By Proposition \ref{pro32}, we have 
	$$S\cap B_\rho(x)\su B_{ \rho^2/5}(L_x^{1'})$$ for some $k$-dimensional affine subspace $L_x^{1'}$.
	Write as above
	$$S\cap B_\rho(x)\su\bigcup_{y\in \3^2_x}B_{\rho}(y)=\bigcup_{y\in \3^2_{x,b}}B_{\rho^2}(y)\cup\bigcup_{y\in \3^2_{x,g}}B_{\rho^2}(y)\equiv B_{\rho^2}(\3^2_{x,b})\bigcup B_{\rho^2}(\3^2_{x,g}),$$
	with $\3^2_x=\3^2_{x,b}\cup\3^2_{x,g}\su L_x^{1'}\cap B_\rho(x)$ so that if $y\neq z$ and $y,z\in\3^2_x$, then $B_{\rho^2/5}(y)\cap B_{\rho^2/5}(z)=\emptyset.$
	Moreover,
	\begin{align*}
		\3^2_{x,g}:=&\{y\in\3^2_{x}:\ F_y^2 \su B_{ \rho^3/5}(L_y^2) \
		\ \ \mbox{\rm for some}\ (k-1)\mbox{\rm -dimensional affine subspace}\ L_y^2\}
	\end{align*}
	with
	\begin{align*}
		F_y^2=\{z\in S\cap B_{2\rho}(y):\Phi_u(z, \rho^3/10)\geq E-\de\}
	\end{align*}
	and
	\begin{align*}
		\3^2_{x,b}:=\{y\in\3^2_x:\ F_y^2\ (\rho^3/10)\mbox{-}\mbox{\rm effectively spans a}\ k\mbox{\rm -dimensional affine subspace}\ L_y^{2'}\}.
	\end{align*}
	Set $\3_b^2=\bigcup\limits_{x\in\3^1_b}\3^2_{x,b}$, $\3^2_g=(\bigcup\limits_{x\in\3^1_b}\3^2_{x,g})\cup\3^1_g$ and $\3^2=\3^2_g\cup\3^2_b$. Then
	$$S\cap B_1(0)\su\bigcup_{x\in \3^2}B_{r^2_{x}}(x)=\bigcup_{x\in \3^2_g}B_{r^2_{x}}(x)\cup\bigcup_{x\in\3^2_b}B_{r^2_{x}}(x)\equiv B_{r^2_{x}}(\3^2_g)\bigcup B_{\rho^2}(\3^2_b),$$
	where $r^2_{x}=\rho$ if $x\in\3^1_g$ and $r^2_{x}=\rho^2$ if $x\in\3^2_g\ba\3^1_g$ or $\3^2_b$.
	Furthermore, with the same reasoning as the initial covering, there holds
	for all $x\in\3^2$ that $\Phi_u(x,r^2_{x})\geq E-\eta$ and
	for all $s\in[r^2_{x},1]$, $u$ is not $(k+1,\ep/2)$-symmetric on $B_{s}(x)$, that is, $x\in S^k_{\ep/2,r^2_{x}}.$
	
	Repeating this procedure, we then build a covering of the form
	\[
	S\cap B_1(0)\su \bigcup_{x\in \3^j}B_{r^j_{x}}(x)=\bigcup_{x\in \3^j_g}B_{r^j_{x}}(x)\cup\bigcup_{\3^j_b}B_{r^j_{x}}(x)\equiv B_{r^j_{x}}(\3^j_g)\bigcup B_{r^j_{x}}(\3^j_b),
	\]
	with $\3_b^j=\bigcup\limits_{x\in\3^{j-1}_b}\3^j_{x,b}$, $\3^j_g=\left(\bigcup\limits_{x\in\3^{j-1}_b}\3^j_{x,g}\right)\cup\3^{j-1}_g$ and $\3^j=\3^j_b\cup\3^j_g,j\geq2.$ Notice that
	\[
	r^j_{x}=\begin{cases}
		\rho^i, & \text{ if } x\in\3^i_g\ba\3^{i-1}_g,i=2,\cdots,j,\\
		\rho^j, & \text{ if } x\in\3^j_b.
	\end{cases}
	\]
	Moreover, $\3^j_{g}$ is defined to be
	\begin{align*}
		\{x\in\3^j:\ r^j_{x}\geq\rho^j\ \mbox{\rm and}\ F_x^j\su B_{ \rho r^j_{x}/5}(L_x^j) \
		\mbox{\rm for some}\ (k-1)\mbox{\rm -dimensional affine subspace}\ L_x^j\}
	\end{align*}
	with
	\begin{align*}
		F_x^j:=\{y\in S\cap B_{2r^j_{x}}(x):\Phi_u(y, \rho r^j_{x}/10)\geq E-\de\}
	\end{align*}
	and
	\begin{align*}
		\3^j_{b}:=\{x\in\3^j:& r^j_{x}=\rho^j \ \mbox{\rm and}
		\ F_x^j\ \frac{\rho r^j_{x}}{10}\mbox{-}\mbox{\rm effectively spans a}\ k\mbox{\rm -dimensional affine subspace}\ L_x^{j'}\}.
	\end{align*}
	Furthermore,  for all $x\neq y\in\3^j\backslash \3_g^{j-1}$, we have $B_{r_x^j/5}(x)\cap B_{r_x^j/5}(y)=\emptyset$ and
	\begin{itemize}
		\item[(P1)] For all $x\in\3^j,$  $\Phi_u(x,r^j_{x})\geq E-\eta.$ 
		\item[(P2)] For all $x\in\3^j$ and for all $s\in[r^j_{x},1]$, $u$ is not $(k+1,\ep)$-symmetric on $B_{s}(x)$, that is, $x\in S^k_{\ep/2,r^j_{x}}.$
	\end{itemize}

	Taking $j=l$ and $\rho^{l}=r$ in \eqref{332} and $\3=\3^{l}$. Then the covering part of Lemma \ref{oooa} is almost complete, except that balls with centers in $\3$ does not necesssarily satisfy the disjointness condition $B_{r_x^i/5}(x)\cap B_{r_y^j/5}(y)=\emptyset$ for all $x,y\in \3$. However, we may do a further covering for the collection of balls $\{B_{r_x}(x)\}$ to select a subcollection of disjoint balls $\{B_{r_x}(x)\}$ such that $\{B_{5r_x}(x)\}$ satisfies all the covering requirements of the lemma. Thus, by relabeling these balls if necesary, we may assume the collection $\{B_{r_x}(x)\}_{x\in \3}$ satisfies the covering part of Lemma \ref{oooa}. 
	\medskip 
	
	\textbf{Step 2. (Reifenberg estimates)}
	In order to prove the volume estimate in \eqref{pla}, we define a measure
	\[
	\mu:=\omega_k\sum_{x\in\3}r_x^k\de_x.
	\]
	and  measures
	\[
	\mu_t:=\omega_k\sum_{x\in\3_t}r_x^k\de_x,\qquad \text{with}\quad \3_t=\{x\in\3:r_x\leq t\},
	\]
	for all $t\in(0,1]$.

	Set $r_j=2^jr,j=0,1,\cdots,al-3$. Then $r_{al-3}=1/8$ by \eqref{332}.
	Since $B_{r/5}(x)\cap B_{r/5}(y)=\emptyset$ for all $x,y\in\3_r,$ it follows easily that $\mu_r(B_r(x))\leq c(m)r^k.$
	Assume now for all  $x\in B_3(0)$ and $s\geq r,$ we have
	\begin{align}\label{1300}
		\mu_{r_j}(B_{r_j}(x)):=\left(\omega_k\sum_{x\in\3,r_x\leq r_j}r_x^k\de_x\right)\left(B_{r_j}(x)\right)\leq C_R(m)r_j^k,
	\end{align}
	where $C_R(m)$ is the constant in Theorem \ref{rre}.  We next show \eqref{1300} is true for ${j+1}$ and hence it holds for all $j\leq al-3$ by induction.
	
	
	We first show that \eqref{1300} holds with constant $C_1(m)=c(m)C_R(m)$. To this end, we write
	\[
	\mu_{r_{j+1}}=\mu_{r_j}+\widetilde{\mu}_{r_{j+1}}:=\sum_{x\in\3_{{r_j}}}\omega_kr_x^k\de_x+\sum_{x\in\3,r_x\in({r_j},{r_{j+1}}]}\omega_kr_x^k\de_x.
	\]
	Take a covering of $B_{r_{j+1}}(x)$ by $M$ balls $\{B_{{r_j}}(y_i)\}$, $M\leq c(m)$, such that $\{B_{{r_j}/5}(y_i)\}$ are disjoint. Then by induction we have
	\[
	\mu_{{r_j}}(B_{r_{j+1}}(x))\leq\sum_{j=1}^{M}\mu_{{r_j}}(B_{{r_j}}(y_i))\leq c(m)C_R(m)r_j^k.
	\]
	By definition of $\widetilde{\mu}_{r_{j+1}}$ and pairwise disjointness of $\{B_{r_x/5}(x)\}$, we have
	\[
	\widetilde{\mu}_{r_{j+1}}(B_{r_{j+1}}(x))\leq c(m)r_{j+1}^k.
	\]
	Thus for all $x\in B_1(0)$, there holds
	\begin{align}\label{13200}
		\mu_{r_{j+1}}(B_{r_{j+1}}(x))\leq c(m)C_R(m)r_{j+1}^k.
	\end{align}
	
	For a fixed ball $B_{r_{j+1}}(x_0)$, set
	\[
	{\mu_{j+1}}=\mu_{r_{j+1}}|_{B_{{r_{j+1}}}(x_0)},
	\]
	and we claim that for all $z\in \supp(\mu_{j+1})$
	\begin{align}\label{137}
		\beta_2=\beta^k_{2,{\mu_{j+1}}}(z,s)^2\leq C_1s^{-k}\int_{B_s(z)}{W}_s(y)d{\mu_{j+1}}(y).
	\end{align}
	
	Note that \eqref{137} trivally holds if $0<s\leq r_z/5$ and we next consider the case $s\geq r_z/5$.  Since  $\Phi_u(z,10s)\leq E$, by (P1) we have 
	$$\Phi_u(z,10s)-\Phi_u(z,5s)\leq\eta\quad \text{ for all } z\in \mbox{supp}(\mu) \text{ and all }s\in[ r_{z}/5,1/10].$$
	Given $\delta_8(m,N,\Lambda,\rho,\ep)$ as  in Theorem \ref{ji}, it follows from Proposition \ref{pro1} that there exists $\eta_0=\eta_0(m,N,\Lambda,\rho,\eta)$ $=\de_1(m,N,\Lambda,\rho,\ep,\de_8)>0$ such that if $\eta\leq\eta_0$, then $u$ is $(0,\de_8)$-symmetric on $B_{10s}(x)$. By (P2), $u$ is not $(k+1,\ep/2)$-symmetric on  $B_{10s}(x)$. Now, \eqref{137} follows from Theorem \ref{ji}. 
	
	We may assume without loss of generality $W_s(x)=0$ if $0<s\leq r_x/5$. Note that, by the induction assumption and \eqref{13200}, for all $j\leq al-3$ and $s\in(0,r_{j+1}],$ and $z\in B_1(0),$ we have
	\begin{equation}\label{eq:rough growth for mu s}
		\mu_s(B_s(z))\leq c(m)C_R(m)s^k.
	\end{equation}
	We claim that for any $r<s\leq r_{j+1}$, we have
	\begin{equation}\label{eq:improved growth for mu j}
	 \mu_{r_{j+1}}(B_s(z))\leq c(m)C_R(m)5^ks^k.
	\end{equation}
	Indeed, if $y\in B_s(z)\cap \supp(\mu)$, then $\frac{r_y}{5}\leq |y-z|\leq s$ and so $y\in \3_{5s}$, which implies $B_s(z)\cap \supp(\mu)\subset \3_{5s}$. Since $r\leq 5s\leq 5r_{j+1}$, we have
	\[
	\mu_{j+1}(B_s(z))\leq \mu_{5s}(B_s(z))\leq \mu_{5s}(B_{5s}(z))\leq c(m)C_R(m)5^ks^k.
	\]
	Integrating \eqref{137} for $s\leq r\leq r_{j+1}$ and $y\in B_{r_{j+1}}(x_0)$ leads to
	\[
	\begin{aligned}
	\int_{B_r(y)} \beta^k_{2,\mu_{j+1}}(z,s)^2d\mu_{j+1}(z)&\leq C_1s^{-k}\int_{B_r(y)}\left[\int_{B_s(z)}{W}_s(x)d\mu_{j+1}(x)\right]d\mu_{j+1}(z)\\
	&\leq C_1s^{-k}\int_{B_r(y)}\int_{B_{2r}(y)}\chi_{B_s(z)}(x)W_s(x)d\mu_{j+1}(x)d\mu_{j+1}(z)\\
	&=C_1s^{-k}\int_{B_{2r}(y)}\mu_{j+1}(B_s(z))d\mu_{j+1}(x)\\
	&\stackrel{\eqref{eq:improved growth for mu j}}{\leq}C_1c(m)C_R(m)\int_{B_{2r}(y)}W_s(x)d\mu_{j+1}(x).
    \end{aligned}
	\]
	A further integration in $s$ from 0 to $r$ leads to
	\[
	\int_{B_r(y)}\int_0^r \beta^k_{2,\mu_{j+1}}(z,s)^2\frac{ds}{s}d\mu_{j+1}(z)\leq c(m)C_1C_R\int_{B_{2r}(y)}\int^r_0{W}_s(z)\frac{ds}{s}d{\mu_{j+1}}(z).
	\]
	Observe that for all  $x\in \mbox{supp}(\mu)$ and $r\leq r_{j+1}\leq1/8,$ we have
	\[
	\int_0^r{W}_s(x)\frac{ds}{s}=\int_{r_x/5}^r{W}_s(x)\frac{ds}{s}\leq\int_{r_x/5}^{1/10}{W}_s(z)\frac{ds}{s}\leq c[\Phi_u(x,1)-\Phi_u(x,r_x/5)]\leq c\eta.
	\]
	Using again the induction hypothesis and \eqref{13200}, we arrive at
	\begin{align}\label{1380}
		\int_{B_r(y)}\left(\int_0^r\beta^k_{2,\mu_{j+1}}(z,s)^2\frac{ds}{s}\right)d\mu_{j+1}(z)\leq cc(m)C_1C_R^2\eta r_{j+1}^k\leq cc(m)C_1C_R^22^{al-3}\eta r^k
	\end{align}
	for all $y\in B_{r_{j+1}}(x)$ and $2^{3-al}r_{j+1}\leq r\leq r_{j+1}.$
	
	Our desired estimate \eqref{137} follows by choosing $\eta $ small enough such that
	$$\eta\leq 2^{3-al}\frac{\de_7^2}{cc(m)C_1C_R^2},$$
	and then applying Theorem \ref{rre} to $\mu_{j+1}$.
\end{proof}


\subsection{The second covering lemma}
The covering obtained in the previous lemma is not sufficient to show the main theorem and so we refine it by deducing the following iterating lemma.
\begin{lemma}[Covering Lemma II]\label{ti}
 Let $u\in H_\Lambda^{2}(B_3,N)$ be a stationary biharmonic map.
  Fix any $\ep>0$ and $0<r<R\leq1$. There exist $\de=\de(m,N,\La,\ep)>0$ and $C_2(m)$ such that  for any subset
 ${S}\su S^k_{\ep,\de r}$,  there exists a finite covering of $S\cap B_R(0)$ by
 \[
 \begin{aligned}
S\cap B_R(0)\su \bigcup_{x\in\3}B_{r_x}(x), \quad\text{with}\quad \ r_x\geq r\quad\text{and}\quad\sum_{x\in\3}r_x^k\leq C_2(m)R^k.
 \end{aligned}
 \]
 Moreover, for each $x\in\3,$ one of the following conditions  is satisfied:
 \begin{itemize}
 	\item[i)] $r_x=r;$
 	\item[ii)] we have the following uniform energy drop property:
 	\begin{equation}\label{110}
 		\sup_{y\in B_{2r_x}(x)\cap S}\Phi_u(y,r_x)\leq E-\de,
 	\end{equation} 
 	where $$E:=\sup_{x\in B_{2R}(0)\cap{S}}\Phi_u(x,R).$$
 \end{itemize}
\end{lemma}

For simplicity we will assume that $R=1$. We will use an induction
argument to refine the first Covering Lemma so as to deduce the desired
covering. The proof is divided into several steps. We will use superscripts
$f,b$ to indicate \emph{final} and \emph{bad} balls respectively. We will also call
the above properties i) and ii) as stopping conditions in our induction
argument below. Note that even though we have the energy drop property,
we can not apply this lemma on $B_{2r_{x}}(x)$ directly due to the
different scale. 
\begin{proof}
\textbf{Step 1. (Initial covering)} By the previous lemma, we have a covering of $S\cap B_{1}(0)$ given by
\[
S\cap B_{1}(0)\subset\bigcup_{x\in{\cal C}_{r}^{0}}B_{r}(x)\cup\bigcup_{x\in{\cal C}_{+}^{0}}B_{r_{x}}(x),
\]
where
\[
{\cal C}_{r}^{0}=\{x\in{\cal C}:r_{x}=r\}\quad\text{and}\quad{\cal C}_{+}^{0}=\{x\in{\cal C}:r_{x}>r\},
\]
and
\begin{equation}\label{eq: 6.67}
\sum_{x\in{\cal C}_{r}^{0}\cup{\cal C}_{+}^{0}}\om_{k}r_{x}^{k}\le C_1(m)=:C_{V}(m).
\end{equation}
Moreover, for each $x\in{\cal C}_{+}^{0}$, the set
\[
F_{x}=\left\{ y\in S\cap B_{2r_{x}}(x):\Phi(y, \rho r_{x}/10)>E-\de\right\} 
\]
is contained in a small neighborhood of a ($k-1$)-dimensional affine subspace. To proceed, we only need to refine the part ${\cal C}_{+}^{0}$. 
\medskip 

\textbf{Step 2. (Recovering of bad balls: the first step)} 

Let $x\in{\cal C}_{+}^{0}$. If $\rho r_{x}=r$, then we just cover $S\cap B_{r_{x}}(x)$
by a family of balls $\{B_{r_{y}}(y)\}_{y\in{\cal C}_{x}^{(1,r)}}$
with $\{B_{\rho r_{x}/2}(y)\}_{y\in{\cal C}_{x}^{(1,r)}}$ being pairwise disjoint, where $y\in B_{r_{x}}(x)$, $r_{y}=\rho r_{x}=r$.  A simple volume comparison argument shows that $\sharp\left({\cal C}_{x}^{(1,r)}\right)\le C(m)\rho^{-n}$, which gives
\[
\sum_{y\in{\cal C}_{x}^{(1,r)}}r_{y}^{k}=\sharp\left({\cal C}_{x}^{(1,r)}\right)(\rho r_{x})^{k}\le C(m)\rho^{k-m}r_{x}^{k}=:C_{r}(m,\rho)r_{x}^{k}.
\]
Thus we collect all such points $x$ and set
\[
{\cal C}^{(1,r)}={\cal C}_{r}^{0}\cup\bigcup_{x\in{\cal C}_{+}^{0},\rho r_{x}=r}{\cal C}_{x}^{(1,r)}.
\]
It follows from \eqref{eq: 6.67} that
\[
\sum_{y\in{\cal C}^{(1,r)}}r_{y}^{k}=\left(\sum_{x\in{\cal C}_{r}^{0}}+\sum_{x\in{\cal C}_{+}^{0}}\sum_{y\in{\cal C}_{x}^{(1,r)}}\right)r^{k}\le\sum_{x\in{\cal C}_{r}^{0}}r^{k}+C_{r}(m,\rho)\sum_{x\in{\cal C}_{+}^{0}}r_{x}^{k}\le C_{V}(m)C_{r}(m,\rho).
\]

Next, suppose $\rho r_{x}>r$. Two cases occur.\textbf{
	Case 1}: $F_{x}=\emptyset$. In this case we simply cover $S\cap B_{r_{x}}(x)$
by balls $\{B_{r_{y}}(y)\}_{y\in{\cal C}_{x}^{(1,f)}}$ centered in
$S\cap B_{r_{x}}(x)$ with $r_{y}=\rho r_{x}$, such that $\{B_{\rho r_{x}/2}(y)\}_{y\in{\cal C}_{x}^{(1,f)}}$
are disjoint. In this case, the energy drop property holds and the number $\sharp\left\{ {\cal C}_{x}^{(1,f)}\right\} $
 is bounded from above by a constant $C(n)\rho^{-m}$, so that
\[
\sum_{y\in{\cal C}_{x}^{(1,f)}}r_{y}^{k}=\sharp\left({\cal C}_{x}^{(1,f)}\right)(\rho r_{x})^{k}\le C(m)\rho^{k-m}r_{x}^{k}=C_{f}(m,\rho)r_{x}^{k}.
\]
Note that actually we have $C_{f}(m,\rho)=C_{r}(m,\rho)$. 

Suppose now \textbf{Case 2} occurs: $F_{x}\neq\emptyset$. We call $B_{r_{x}}(x)$ a bad ball. We want to recover $B_{r_{x}}(x)$ using the fact that $F_{x}\subset B_{\rho r_{x}/5}(L_{x}^{k-1})\cap B_{2r_{x}}(x)$.
This will be done as follows. First note that the part away from $B_{\rho r_{x}}(F_{x})$
is good: we have
\[
S\cap B_{r_{x}}(x)\Big\backslash B_{\rho r_{x}}(F_{x})\subset\bigcup_{y\in{\cal C}_{x}^{(1,f)}}B_{r_{y}}(y)\quad\text{with }r_{y}=\rho r_{x},
\]
and $\{B_{\rho r_{x}/2}(y)\}_{{\cal C}_{x}^{(1,f)}}$ being pairwise
disjoint. Since $y$ is away from  $F_{x}$,  the energy drop property holds:
$\Phi_{r_{y}/10}(y)\le E-\de$ for all $y\in{\cal C}_{x}^{(1,f)}$.
Thus these balls will be part of the ``final'' balls in the next
step. Moreover, we have the trivial estimate via a volume comparison
argument:
\[
\sum_{y\in{\cal C}_{x}^{(1,f)}}r_{y}^{k}=\sharp\left\{ {\cal C}_{x}^{(1,f)}\right\} (\rho r_{x})^{k}\le C(m)\rho^{k-n}r_{x}^{k}=C_{f}(m,\rho)r_{x}^{k}.
\]
Then we collect all the ${\cal C}_{x}^{(1,f)}$ from the above to
obtain final balls of the first generation:
\[
{\cal C}^{(1,f)}=\bigcup_{x\in{\cal C}_{+}^{0},\rho r_{x}>r}{\cal C}_{x}^{(1,f)}.
\]
It follows from \eqref{eq: 6.67}  that
\[
\sum_{y\in{\cal C}^{(1,f)}}r_{y}^{k}=\sum_{x\in{\cal C}_{+}^{0}}\sum_{y\in{\cal C}_{x}^{(1,f)}}(\rho r_{x})^{k}\le C_{f}(m,\rho)\sum_{x\in{\cal C}_{+}^{0}}r_{x}^{k}\le C_{V}(m)C_{f}(m,\rho).
\]

For the remaining part of $S\cap B_{r_{x}}(x)$, we simply
cover it by
\[
S\cap B_{r_{x}}(x)\cap B_{\rho r_{x}}(F_{x})\subset\bigcup_{y\in{\cal C}_{x}^{1,b}}B_{r_{y}}(y)\quad\text{with }r_{y}=\rho r_{x},
\]
and $\{B_{\rho r_{x}/2}(y)\}_{{\cal C}_{x}^{(1,b)}}$ being pairwise
disjoint, where ``$b$'' means bad balls. The problem is that the
energy drop condition can not be verified on each ball $B_{r_{y}}(y)$.
However, since $F_{x}\subset B_{\rho r_{x}/5}(L_{x})\cap B_{2r_{x}}(x)$
for some ($k-1$)-dimensional affine subspace $L_{x}$, a volume comparison
argument gives 
\[
\sharp\{{\cal C}_{x}^{(1,b)}\}\le C(m)\rho^{1-k},
\]
This means that there are relatively very few bad balls. Hence
\[
\sum_{y\in{\cal C}_{x}^{(1,b)}}(\rho r_{x})^{k}=(\rho r_{x})^{k}\sharp\left\{ \text{\ensuremath{{\cal C}_{x}^{(1,b)}}}\right\} \le C(m)\rho r_{x}^{k}=C_{b}(m)\rho r_{x}^{k}.
\]
Thus we collect all the bad balls and define
\[
{\cal C}^{(1,b)}=\bigcup_{x\in{\cal C}_{+}^{0},\rho r_{x}>r}{\cal C}_{x}^{(1,b)}.
\]

Now we have the following properties for balls in ${\cal C}^{(1,b)}$
and ${\cal C}^{(1,r)}\cup{\cal C}^{(1,f)}$. For bad balls,  using \eqref{eq: 6.67}, we obtain
\[
\sum_{y\in{\cal C}^{(1,b)}}r_{y}^{k}=\sum_{x\in{\cal C}_{+}^{0}}\sum_{y\in{\cal C}_{x}^{(1,b)}}(\rho r_{x})^{k}\le C_{b}(m)\rho\sum_{x\in{\cal C}_{+}^{0}}r_{x}^{k}\le C_{V}(m)C_{b}(m)\rho.
\]
Hereafter, we choose
\[
0<\rho<\min\left\{ 100^{-1},\frac{1}{2C_{V}(m)C_{b}(m)}\right\} 
\]
such that
\[
\sum_{y\in{\cal C}^{(1,b)}}(\rho r_{x})^{k}<1/2.
\]
Also, by setting
\[
C_{2}=C_{2}(m)=2C_{V}(m)\left(C_{r}(m,\rho)+C_{f}(m,\rho)\right),
\]
we have
\[
\begin{aligned}\sum_{y\in{\cal C}^{(1,r)}\cup{\cal C}^{(1,f)}}r_{y}^{k} & \le(C_{f}(m,\rho)+C_{r}(m,\rho))\sum_{x\in{\cal C}}r_{x}^{k}\le\frac{1}{2}C_{2}(m).\end{aligned}
\]

We also note that for all $y\in{\cal C}^{(1,b)}$, $r<r_{y}\le\rho$. 
\medskip 

\textbf{Step 3. (Induction step)} Now our aim is to derive the following
covering: for each $i\ge1$,
\[
S\cap B_{1}(0)\subset\bigcup_{x\in{\cal C}^{(i,r)}}B_{r}(x)\cup\bigcup_{x\in{\cal C}^{(i,f)}}B_{r_{x}}(x)\cup\bigcup_{x\in{\cal C}^{(i,b)}}B_{r_{x}}(x)
\]
with the following properties hold:

(1) For each $x\in{\cal C}^{(i,r)}$, $r_{x}=r$;

(2) For each $x\in{\cal C}^{(i,f)}$, the energy drop condition holds:
for all $z\in \cal{S}\cap B_{2r_{x}}(x)$ we have $\Phi(z,r_{x}/10)\le E-\de$;
these balls will be called final balls;

(3) For each $x\in{\cal C}^{(i,b)}$, we have $r<r_{x}\le\rho^{i}$.
On these ``bad'' balls, none of the above two stopping conditions is
verified. 

(4) There holds
\begin{equation}
\sum_{y\in{\cal C}^{(i,r)}\cup{\cal C}^{(i,f)}}r_{y}^{k}\le\left(\sum_{j=1}^{i}2^{-j}\right)C_{2}(m)\quad \text{and}\quad \sum_{y\in{\cal C}^{(i,b)}}r_{y}^{k}\le2^{-i}.\label{eq: 6.70}
\end{equation}

We have proved the above induction for $i=1$ in Step 2. Suppose now
it holds for some $i\ge1$. We aim to prove that it holds for
$i+1$. It is clear that we only need to recover these bad balls. 

Suppose now $x\in{\cal C}^{(1,b)}$, i.e., $B_{r_{x}}(x)$ is a bad
ball and $r<r_{x}\le\rho^{i}$. 

\textbf{Case 1: $\rho r_{x}=r$.} Then we cover it simply by balls
of radius $\rho r_{x}$ as in step 2, and obtain a covering of balls
whose centers lie in ${\cal C}_{x}^{(i+1,r)}$, together with the estimate
\[
\sum_{y\in{\cal C}_{x}^{(i+1,r)}.}(\rho r_{x})^{k}\le C(m)\rho^{k-n}r_{x}^{k}=C_{r}(m,\rho)r_{x}^{k}.
\]

\textbf{Case 2: $\rho r_{x}>r$.} In this case we apply the previous lemma 
to get a covering $\{B_{r_{y}}(y)\}_{y\in{\cal C}^{x}}$ such that
\[
S\cap B_{r_{x}}(x)\subset\bigcup_{y\in{\cal C}_{r}^{x}}B_{r}(y)\cup\bigcup_{y\in{\cal C}_{+}^{x}}B_{r_{y}}(y)\qquad\text{with }r_{y}\ge r
\]
and
\[
\sum_{y\in{\cal C}_{r}^{x}\cup{\cal C}_{+}^{x}}\om_{k}r_{y}^{k}\le C_{V}(m)r_{x}^{k}.
\]
Moreover, for each $y\in{\cal C}_{+}^{x}$, there is a ($k-1$)-dimensional
affine subspace $L_{y}$ such that
\[
F_{y}\equiv\left\{ z\in S\cap B_{2r_{y}}(y):\Phi_{\rho r_{y}/10}(z)\ge E-\de\right\} \subset B_{\rho r_{y}/5}(L_{y})\cap B_{2r_{y}}(y).
\]
Of course we will reserve ${\cal C}_{r}^{x}$ as part of ${\cal C}^{(i+1,r)}$.
Thus below we assume that $y\in{\cal C}_{r}^{x}$. 

\textbf{Case 2.1: $\rho r_{y}=r$.} Then as that of Case 1, we get
a simple covering of at most $C(m)\rho^{-m}$ $r$-balls $\{B_{r}(z)\}_{z\in{\cal C}_{y}^{(i+1,r)}}$
of $B_{r_{y}}(y)$. Then we define
\[
{\cal C}^{(i+1,r)}={\cal C}^{(i,r)}\cup\bigcup_{x\in{\cal C}^{(i,b)},\rho r_{x}=r}{\cal C}_{x}^{(i+1,r)}\cup\bigcup_{x\in{\cal C}^{(i,b)},\rho r_{x}>r}\left({\cal C}_{r}^{x}\cup\bigcup_{y\in{\cal C}_{+}^{x},\rho r_{y}=r}{\cal C}_{y}^{(i+1,r)}\right).
\]
This shows that how much more is ${\cal C}^{(i+1,r)}$ than that of
${\cal C}^{(i,r)}$. We get
\[
\begin{aligned}\sum_{z\in{\cal C}^{(i+1,r)}}r_{z}^{k} & =\sum_{z\in{\cal C}^{(i,r)}}r_{z}^{k}+\sum_{x\in{\cal C}^{(i,b)},\rho r_{x}=r}\sum_{z\in{\cal C}_{x}^{(i+1,r)}}r_{z}^{k}\\
& \quad+\sum_{x\in{\cal C}^{(i,b)},\rho r_{x}>r}\left(\sum_{z\in{\cal C}_{r}^{x}}r_{z}^{k}+\sum_{y\in{\cal C}_{+}^{x}}\sum_{z\in{\cal C}_{y}^{(i+1,r)}}r_{z}^{k}\right)\\
& \le\sum_{z\in{\cal C}^{(i,r)}}r_{z}^{k}+C_{r}(m,\rho)\left(\sum_{x\in{\cal C}^{(i,b)},\rho r_{x}=r}r_{x}^{k}+\sum_{x\in{\cal C}^{(i,b)},\rho r_{x}>r}\sum_{y\in{\cal C}^{x}}r_{y}^{k}\right)\\
& \le \sum_{z\in{\cal C}^{(i,r)}}r_{z}^{k}+C_{r}(m,\rho)C_{1}(m)\sum_{x\in{\cal C}^{(i,b)}}r_{x}^{k}\\
& \le\sum_{z\in{\cal C}^{(i,r)}}r_{z}^{k}+2^{-i}C_{r}(m,\rho)C_{1}(m).
\end{aligned}
\]

\textbf{Case 2.2: $\rho r_{y}>r$.} Then whether $F_{y}$ is empty
or not, we cover $B_{r_{y}}(y)$ by
\[
\begin{aligned} & S\cap B_{r_{y}}(y)\Big\backslash B_{\rho r_{y}}(F_{y})\subset\bigcup_{z\in{\cal C}_{y}^{(i+1,f)}}B_{r_{z}}(z),\\
& S\cap B_{r_{y}}(y)\cap B_{\rho r_{y}}(F_{y})\subset\bigcup_{z\in{\cal C}_{y}^{(i+1,b)}}B_{r_{z}}(z),
\end{aligned}
\quad\text{with }r_{z}=\rho r_{y},
\]
and each ball $B_{r_{z}}(z)$ with $z\in{\cal C}_{y}^{(i+1,f)}$ satisfies
the energy drop condition. Moreover, the following estimates hold:
\[
\sum_{z\in{\cal C}_{y}^{(i+1,f)}}r_{z}^{k}\le C_{f}(m,\rho)r_{y}^{k}\quad\text{and}\quad \sum_{z\in{\cal C}_{y}^{(i+1,b)}}r_{z}^{k}\le C_{b}(m)\rho r_{y}^{k}.
\]
Note that in this case no $r$-balls occur anymore. Thus we set
\[
{\cal C}^{(i+1,f)}={\cal C}^{(i,f)}\cup\bigcup_{x\in{\cal C}^{(i,b)}}\bigcup_{y\in{\cal C}_{+}^{x}}{\cal C}_{y}^{(i+1,f)},\quad{\cal C}^{(i+1,b)}=\bigcup_{x\in{\cal C}^{(i,b)}}\bigcup_{y\in{\cal C}_{+}^{x}}{\cal C}_{y}^{(i+1,b)}.
\]
Then we have
\[
\begin{aligned}\sum_{z\in{\cal C}^{(i+1,f)}}r_{z}^{k} & =\sum_{z\in{\cal C}^{(i,f)}}r_{z}^{k}+\sum_{x\in{\cal C}^{(i,b)}}\sum_{y\in{\cal C}_{+}^{x}}\sum_{z\in{\cal C}_{y}^{(i+1,f)}}r_{z}^{k}\\
& \le\sum_{z\in{\cal C}^{(i,f)}}r_{z}^{k}+C_{f}(m,\rho)\sum_{x\in{\cal C}^{(i,b)}}\sum_{y\in{\cal C}_{+}^{x}}r_{y}^{k}\\
& \le\sum_{z\in{\cal C}^{(i,f)}}r_{z}^{k}+C_{f}(m,\rho)C_{V}(m)\sum_{x\in{\cal C}^{(i,b)}}r_{x}^{k}\\
& \le\sum_{z\in{\cal C}^{(i,f)}}r_{z}^{k}+2^{-i}C_{f}(m,\rho)C_{V}(m).
\end{aligned}
\]
Recall that $C_{2}(m)=2\left(C_{r}(m,\rho)+C_{f}(m,\rho)\right)C_{V}(m)$.
Hence by the inductive assumption we obtain
\[
\left(\sum_{z\in{\cal C}^{(i+1,r)}}+\sum_{z\in{\cal C}^{(i+1,f)}}\right)r_{z}^{k}\le\left(\sum_{z\in{\cal C}^{(i,r)}}+\sum_{z\in{\cal C}^{(i,f)}}\right)r_{z}^{k}+2^{-i-1}C_{2}(m)\le C_{2}(m)\sum_{j=1}^{i+1}2^{-j},
\]
and also by the choice of $\rho$, we obtain
\[
\begin{aligned}\sum_{z\in{\cal C}^{(i+1,b)}}r_{z}^{k}= & \sum_{x\in{\cal C}^{(i,b)}}\sum_{y\in{\cal C}_{+}^{x}}\sum_{z\in{\cal C}_{y}^{(i+1,b)}}r_{z}^{k}\le\sum_{x\in{\cal C}^{(i,b)}}\sum_{y\in{\cal C}_{+}^{x}}C_{b}(m)\rho r_{y}^{k}\\
& \le\sum_{x\in{\cal C}^{(i,b)}}C_{b}(m)\rho C_{1}(m)r_{y}^{k}\le2^{-1}\sum_{x\in{\cal C}^{(i,b)}}r_{y}^{k}\le2^{-i-1}.
\end{aligned}
\]
This proves \eqref{eq: 6.70} for $i+1$. 

To finish the proof, we only need to note that for every $z\in{\cal C}^{(i+1,b)}$,
we have $r_{z}=\rho r_{y}\le\rho r_{x}$ for some $x\in{\cal C}^{(i,b)}$,
and hence $r_{z}\le\rho^{i+1}$ since $r_{x}\le\rho^{i}$ by inductive
assumption. The proof of Step 3 is complete.

Since $r=\rho^{\bar{j}}$, the above procedure will stop at $i=\bar{j}-1$.  
\medskip 

\textbf{Step 4. (Final refinement) } Now we have obtained a covering of $S\cap B_R(0)$ satisfying 
\[
\begin{aligned}
S\cap B_R(0)\su \bigcup_{x\in\3}B_{r_x}(x)=\bigcup_{x\in\3_{r}}B_{r_x}(x)\cup\bigcup_{x\in\3_{+}}B_{r_x}(x),\quad  \sum_{x\in\3}r_x^k\leq C_2(m),
\end{aligned}
\]
where $\3_{r}$ consists of centers of all balls $B_{r_x}(x)$ such that $r_x=r$ and  $\3_{+}$ consists of centers of all balls $B_{r_x}(x)$ such that $r_x> r$ and
\[
\begin{aligned}
\sup_{y\in S\cap B_{2r_x}(x)}\Phi_u({y,r_x/10})\leq E-\de.
\end{aligned}
\]
For each $x\in\3_{+}$, we can simply cover  $S\cap B_{r_x}(x)$ by a family of balls $\{B_{\rho r_x}(y)\}_{y\in \3_x}$  such that $y\in S\cap B_{r_x}(x) $,  $\{B_{\rho r_x/5}(y)\}_{y\in \3_x}$ are pairwise disjoint and $\sharp\{\3_x\}\leq C(m)\rho^{-m}$. For each $y\in \3_x$, set $r_y=\rho r_x$. Then we have $B_{2r_y}(y)\subset  B_{2r_x}(x)$ and $r_y\le r_x/10$, which implies by monotonicity that
\[
\sup_{z\in S\cap B_{2r_y}(y)}\Phi_u(z,r_y)\leq \sup_{z\in S\cap B_{2r_x}(x)}\Phi_u(z,r_x/10)\leq E-\delta. 
\]
Moreover, 
\[
\sum_{x\in \3_+}\sum_{y\in \3_x}r_y^k\leq C(m)\rho^{k-m}\sum_{x\in \3_+}r_x^k\leq C(m)\rho^{k-m}C_2(m).
\]
This completes the proof of Lemma \ref{ti} upon taking $C_2'(m)=2C(m)\rho^{k-m}C_2(m)$.
\end{proof}

Now we can complete the proof of Lemma \ref{lemma:main covering lemma}.
\begin{proof}[Proof of Lemma \ref{lemma:main covering lemma}]
Note that the energy $E$ defined as in Lemma \ref{ti} satisfies $E\leq C\Lambda$.  So iterating Lemma \ref{ti} by at most $i=[\de^{-1}E]+1$ times, we could obtain a covering $\{B_{r_x}(x)\}_{x\in \3^i}$ of $S\cap B_R(0)$ such that $r_x\le r$ and 
\[
\sum_{x\in \3^i}r_x^k\leq C_F(m)^iR^k.
\] 
We may further assume that $x\in S\cap B_R(0)$ by considering the larger covering $\{B_{2r_x}(x)\}_{x\in \tilde{\cal C}^i}$. Then 
\[ S\cap B_R(0) \subset \bigcup_{x\in \tilde{\cal C}^i} B_{2 r_x}(x)\qquad \text{ and } \quad  \sum_{x\in \cal{C}}(2 r_x)^k\leq 2^kC_F(m)^iR^k.\]
 Finally, since $\sup\limits_{x\in \tilde{\cal C}^i}r_x\le r$, we can use Vitali's covering lemma  to select a family of disjoint balls $\{B_{2r_x}(x)\}_{x\in \cal{C}}$ from that of $\tilde{\cal C}^i$ so that
 \[ S\cap B_R(0) \subset \bigcup_{x\in \cal{C}} B_{10r_x}(x)\qquad \text{ and } \quad  \sum_{x\in \cal{C}}(10r_x)^k\leq 10^kC_F(m)^iR^k.\]
  The proof is complete upon taking $C(m)=10^kC_F(m)^i$ and relabelling the balls.
\end{proof}

\section{Proofs of main theorems}
 \subsection{Proof of Theorem \ref{thm80}}
This theorem follows from the covering lemmas and the rectifiable Reifenberg theorems mentioned in the previous sections.
\begin{proof}[Proof of Theorem \ref{thm80}]

By Lemma \ref{lemma:main covering lemma}, for all $0<r<\delta$,
\begin{align}\label{z3}
	\mbox{Vol}(T_r(S^k_{\ep,r}(u))\cap B_1(0))\le \mbox{Vol}\Big(T_r\Big(S^k_{\ep,\de }(u)\cap B_1(0)\Big)\Big) \leq C'_\ep r^{m-k}
\end{align}
with $C'_{\ep}=C'_\ep(m,N,\Lambda,\ep)$.   This proves the desired volume estimate \eqref{J} for all $0<r<\de$. In general, since $\de$ also depends on $m,N,\Lambda,\ep$, for any $\de\le r<1$ we have 
\[\label{z3}
\mbox{Vol}\Big(T_r\Big(S^k_{\ep,r}(u)\cap B_1(0)\Big)\Big) \le \mbox{Vol}(T_1( B_1(0)) \leq C_m \le C_m \left(\frac{r}{\de} \right)^{m-k}.\]
This completes the proof of  \eqref{J} for all $0<r<1$. The volume estimate \eqref{Jk} follows from  \eqref{J} by noting that $S^k_{\ep}(u)\subset S^k_{\ep, r}(u)$ for any $r>0$.

 In order to prove  the rectifiability of $S^k(u)$, it is sufficient to prove the rectifiability of $S^k_\ep(u)$ for each $\ep>0$, as $S^k(u)=\bigcup\limits_{i\ge 1}S^k_{1/i}(u)$.
 
 By the volume estimate \eqref{Jk}, we have $\HH^k(S^k_\ep(u)\cap B_1(0))\leq C_\ep.$ Applying the same estimates on $B_r(x)$ with $x\in B_1(0)$ and $r\leq 1$ gives the Alhfors upper bound estimate
\begin{equation}\label{jkj}
\HH^k(S^k_\ep(u)\cap B_r(x))\leq C_\ep r^k.
\end{equation}

 Let $S\su S_\ep^k(u)\cap B_1(0)$ be an arbitrary measurable subset with $\HH^k(S)>0$. Set 
 $$g(x,r)=\Phi_u(x,r)-\Phi_u(x,0),\quad x\in B_1(0),\ r\leq 1.$$ Then the dominated convergence theorem implies that for all $\de>0,$ there exists $\overline{r}>0$ such that
 \begin{equation*}\label{139}
\medint_S g(x,10\overline{r})d\HH^k(x)\leq\de^2.
\end{equation*}
So we can find a measurable subset $E\su S$ with $\HH^k(E)\leq\de\HH^k(S)$ and $g(x,10\overline{r})\leq\de$ for all $x\in F:=S\ba E$. Cover $F$ by a finite number of balls $B_{\overline{r}}(x_i)$ centered on $F.$ Rescaling if necesssary, we may assume that $B_{\overline{r}}(x_i)=B_1(0)$. Then $g(x,10)\leq\de$ for $x\in F.$ Similar to \eqref{137}, choose $\de$ sufficiently small such that $u$ is $(0,\de_9)$ symmetric in $B_{10}$. Theorem \ref{ji} implies 
\begin{equation*}
\beta_{2,\HH^k|_F}(z,s)^2\leq C_1s^{-k}\int_{B_s(z)}W_s(t)d\HH^k|_F(t) \quad \text{for all }z\in F,\ s\leq1.
\end{equation*}
Integrating the previous estimate with respect to $z$ and using \eqref{jkj} yield that for all  $x\in B_1(0)$ and $s\leq r\leq1,$ there holds
\begin{align*}
\int_{B_r(x)}\beta_{2,\HH^k|_F}(z,s)^2d\HH^k|_F(z)&\leq C_1s^{-k}\int_{B_r(x)}\int_{B_s(z)}W_s(t)d\HH^k|_F(t)d\HH^k|_F(z)\non\\
&\leq C_1C_\ep\int_{B_{r+s}(x)}W_s(z)d\HH^k|_F(z).
\end{align*}
Integrating again with respect to $s$, similar to \eqref{1380}, we obtain that for all $x\in B_1(0)$ and $r\leq1,$
\begin{align*}
\int_{B_r(x)}\int_0^r\beta_{2,\HH^k|_F}(z,s)^2\frac{ds}{s}d\HH^k|_F(z)&\leq C_1C_\ep\int_{B_{2r}(x)}[\Phi_u(z,10r)-\Phi_u(z,0)]d\HH^k|_F(z)\non\\
&\leq c(m)C_1C_\ep^2\de r^k.
\end{align*}
Choosing
$$\de\le\frac{\de_8^2}{c(m)C_1C_\ep^2},$$
 we deduce from Theorem \ref{rree} that $F\cap B_1(0)$ is $k$-rectifiable.
 
Repeating the above argument with $E$ in place of $F$, we could find another measurable set $E_1\subset E$ with $\HH^k(E_1)\leq \delta \HH^k(E)$, and that $F_1:=E\ba E_1$ is $k$-rectifiable. Contiuing this process, we eventually conclude that $S$ is $k$-rectifiable. 

The proofs of final assertions are similar to that of \cite{Naber-V-2017} and are thus omitted.
\end{proof}

\subsection{Proof of Theorem \ref{thm84} and Theorem \ref{thm85}}  The proofs of these two theorems
   are similar to that of \cite{Breiner-Lamm-2015}, which are built on the important compactness results of Scheven \cite{Scheven-2008-ACV}. For the convenience of readers, we sketch their proofs here.  First of all, we prove the following theorem.

\begin{theorem}[{\cite[Theorem  4.9]{Breiner-Lamm-2015}, symmetry implies regularity}] \label{thm: CN-ep-regu}
	There exists $\de(m,K_{N},\La)>0$ such that
	\[
	r_{u}(0)\ge1/2
	\]
 for any biharmonic map $u\in H^2_\Lambda(B_1,N)$,	provided one of the following conditions is satisfied:
 \begin{itemize}
 	\item[(1)]  $u$ is an $(m-4,\de)$-symmetric minimizing biharmonic map;
 	
 	\item[(2)] The target manifold $N$ does not admit any nonconstant smooth 0-homogeneous stationary biharmonic map from
 	$\R^5\backslash\{0\}$ to $N$, and $u$ is an $(m-5,\de)$-symmetric stationary biharmonic map;
 	
 	\item[(3)] For some $k\ge 4$, $N$ does not admit any nonconstant smooth 0-homogeneous minimizing biharmonic map from 
 	$\R^{l+1}\backslash\{0\}$ to $N$ for all $4\le l\le k$, and $u$ is an $(m-k-1,\de)$-symmetric minimizing biharmonic map;
 	
 	\item[(4)] For some $k\ge 4$,  $N$ does not admit any nonconstant smooth 0-homogeneous sationary biharmonic map from 
 	$\R^{l+1}\backslash\{0\}$ to $N$ for all $4\le l\le k$, and $u$ is an $(m-k-1,\de)$-symmetric stationary  biharmonic map.
 \end{itemize}
%
%
%
%
	\end{theorem}

Note the this theorem is different from the usual partial regularity 
theory, since it only assumes that the map satisfies certain quantitative $\ep$-symmetry.
Condition (2) can be unified by Condition (4). But we keep it as a prototype of Condition (4). The proof is divided into the following three lemmas, which are of independent interest. The first one is a compactness result, which are collected from Theorem 1.5, Thoerem 1.6 and Proposition 1.7 of Scheven \cite{Scheven-2008-ACV} and Lemma 4.7 of Breiner and Lamm \cite{Breiner-Lamm-2015}.
 \begin{lemma}[Compactness]\label{lem: compactness}  
 	Let $\{u_i\}_{i\ge 1}\subset H^2_{\La}(B_{1},N)$ be a sequence of uniformly bounded stationary biharmonic map, and $u_i\wto u$ in $H^2$. Then $u$ is a weakly biharmonic map. Moreover, there holds
 	\[u_i\to u\qquad \text{strongly in }H^2,\] provided one of the following conditions is satisfied:
 	\begin{itemize}
 		\item[(i)] $\{u_i\}_{i\geq 1}$ is a sequence of minimizing biharmonic maps. In this case,  $u$ is also a minimizing biharmonic map, provided $m=5$ and $u$ is 0-homogeneous, or $N$ satisfies condition {\upshape(3)} of Theorem \ref{thm: CN-ep-regu} and $u$ is $(m-l)$-symmetric for any $4\le l\le k+2$;
 		
 		\item[(ii)] The target manifold $N$ does not admit any nonconstant smooth 0-homogeneous stationary biharmonic map from 
 		$\R^5\backslash\{0\}$ to $N$. In this case, $u$ is also a stationary biharmonic map.
 	\end{itemize}
%
%
\end{lemma}

We remark that  under Condition (i) it is still \emph{unknown} whether the limit map $u$ is  minimizing biharmonic  without any further assumptions.  Due to this lemma, it is  natural to  ask   whether a bounded sequence of weakly biharmonic maps would converge weakly (up to a subsequence) to a weakly  biharmonic maps.  

The compactness result implies the following regularity result. 	

\begin{lemma} [{$(m,\ep)$-regularity theory}]\label{lemma:new epsilon regularity} 
	There exists $\ep>0$ depending only on $m,\La,N$ such that
    $$r_{u}(0)\ge1/2,$$ provided one of the following conditions is satisfied:
    \begin{itemize}
    	\item[(1)] $u\in H^2_{\La}(B_{1},N)$ is an $(m,\ep)$-symmetric minimizing biharmonic map;
    	
    	\item[(2)] $u\in H^2_{\La}(B_{1},N)$ is an $(m,\ep)$-symmetric stationary biharmonic map and the  target manifold $N$ does not admit any nonconstant smooth 0-homogeneous stationary biharmonic map from $\R^5\backslash\{0\}$ to $N$.
    \end{itemize}
	
%
	\end{lemma}

 Note also that under Condition (2), the $(m,\ep)$-regularity theory for stationary biharmonic maps  is an improvement of the usual partial regularity Theorem \ref{thm3}).
\begin{proof}
	Assume, on the contrary, that there is a sequence of minimizing (stationary
	resp.) biharmonic maps $u_{i}\in H_{\La}^{2}(B_{1},N)$ satisfying
	\[
	\medint_{B_{1}}|u_{i}-c_{i}|^{2}dx<1/i,
	\]
	but $u_{i}\not\in C^{\wq}(B_{1/2})$. Up to a further subsequence, we may assume that $u_{i}\wto u_{\wq}$
	in $H^{2}(B_{1})$ and $u_{i}\to u_{\wq}$ in $L^{2}(B_{1})$. Sending $i\to\wq$, we find that $u_{\wq}\equiv$ constant. Moreover, under either one of the two conditions, 
	 the compactness  Lemma \ref{lem: compactness} implies the strong convergence $u_{i}\to u_{\wq}$ in $H^{2}$. As a result,
	\[
	\medint_{B_{1}}|\na^2 u_{i}|^{2}+|\na u_{i}|^{2}\to\medint_{B_{1}}|\na^2 u_{\wq}|^{2}+|\na u_{\wq}|^{2}=0.
	\]
	Then the $\ep$-regularity theory for stationary biharmonic maps implies  $u_{i}\in C^{\wq}(B_{1/2})$
	for $i\gg1$. This is a contradiction.
\end{proof}

The last lemma concerns symmetry improvement.
\begin{lemma}[Symmetry improvement]\label{lemma:symmetry improvement}
	For any $\ep>0$, there exists
	$\de=\de(m,K_{N},\La,\ep)>0$ such that $u$ is $(m,\ep)$-symmetric
	on $B_{1}$, provided $u$ satisfies one of the four conditions in Theorem \ref{thm: CN-ep-regu}.
\end{lemma}
\begin{proof}
		We only prove the result under Condition (1) and (3) of Theorem \ref{thm: CN-ep-regu}, since the proofs for other cases are entirely similar.  We shall use a contradiction argument to prove that $u$ is $(m,\ep)$-symmetric
	on $B_{1}$, that is, there is a constant $c>0$ such that
	\begin{equation}\label{eq: (n,ep)-symmetry}
	\medint_{B_{1}}|u-c|^{2}dx<\ep.
	\end{equation}

	Case I: Condition (1) holds.
	\smallskip 

     Suppose that there is a constant $\ep>0$
	such that, for each $i\ge1$ there is a minimizing biharmonic map $u_{i}\in H_{\La}^{2}(B_{1},N)$
	which is $(m-4,1/i)$-symmetric in $B_{1}$ but not $(m,\ep)$-symmetric
	in $B_{1}$. Then up to a subsequence we have $u_{i}\wto u_{\wq}$
	in $H^{2}(B_{1})$ and $u_{i}\to u_{\wq}$ in $L^{2}(B_{1})$. Moreover,
	$u_{\wq}$ is $(m-4)$-symmetric but not $(m,\ep/2)$-symmetric on
	$B_{1}$.
	
By the compactness  Lemma \ref{lem: compactness}, we know that $u_{\wq}$ is a weakly $(m-4)$-symmetric biharmonic map (we do not need the strong convergence here). Thus $u_{\wq}$ can be viewed as a 0-homogeneous weakly biharmonic map from
	$B^{4}$ to $N$. Hence $u_{\wq}$ is smooth by the regularity theory of 4-dimensional biharmonic maps, see e.g.~\cite{Chang-W-Y-1999,Strzelecki-2003-CV,Wang-2004-CPAM}. But then the 0-homogeneity of $u_{\wq}$ implies that it
	is a constant. This leads to a contradiction, since $u_{\wq}$ is not $(m,\ep/2)$-symmetric
	on $B_{1}$.
	
	\smallskip	
	Case II: Condition (3) holds.
   \smallskip 
   
	Suppose that there is a constant $\ep>0$
	such that, for each $i\ge1$ there is a minimizing harmonic map $u_{i}\in H_{\La}^{2}(B_{1},N)$ which is $(m-k-1,1/i)$-symmetric in $B_{1}$ but not $(m,\ep)$-symmetric in $B_{1}$. Then up to a subsequence we have $u_{i}\wto u_{\wq}$ in $H^{2}(B_{1})$ and $u_{i}\to u_{\wq}$ in $L^{2}(B_{1})$. Moreover, $u_{\wq}\colon B_1\to N$ is $(m-k-1)$-symmetric but not $(m,\ep/2)$-symmetric on $B_{1}$. 
	
	By the compactness  Lemma \ref{lem: compactness}, $u_{\wq}$ is a weakly $(m-k-1)$-symmetric
biharmonic map. Thus $u_{\wq}$ can be viewed as a 0-homogeneous weakly biharmonic map from $B^{k+1}$ to $N$. 	
	To reach a contradiction, we will prove that $u_{\wq}$ is constant. 	To this end, we need to use the assumption on $N$, i.e., $N$ has
	no nonconstant smooth 0-homogeneous minimizing biharmonic map from $\R^{l+1}\backslash\{0\}$ to $N$ for all $4\le l\le k$.  Under this
	assumption, the compactness  Lemma \ref{lem: compactness} implies that $u_{i}\to u_{\wq}$ strongly in $H^{2}$ and $u_{\wq}$
	is a minimizing biharmonic map. Furthermore, by \cite[Theorem 1.1]{Scheven-2008-ACV}, $\dim_{\mathcal{H}}(\sing(u_{\infty}))=0$, and thus $u_{\wq}\in C^\infty(B^{k+1}\backslash\{0\})$ is a 0-homogeneous minimizing biharmonic map, which in turn implies that $u_{\wq}$ is constant by the assumption.  	
\end{proof}

\begin{proof}[Proof of Theorem \ref{thm: CN-ep-regu}]
It follows straightforward from Lemmas \ref{lem: compactness}, \ref{lemma:new epsilon regularity} and \ref{lemma:symmetry improvement}.
\end{proof}

Now we sketch the proof of Theorem \ref{thm84}. 
\begin{proof}[Proof of Theorem \ref{thm84}]
By a scaling argument, 	Theorem \ref{thm: CN-ep-regu} implies
	$$\{x\in B_1(0) : r_u(x)<r\}\su S_{\ep,r}^{m-5}(u).$$
	Thus by Theorem \ref{thm80} there exists $C>0$ such that for each $0<r<1$ we have
	\[
	{\rm Vol}(B_r(\{x\in B_1(0):r_u(x)<r\}))\leq{\rm Vol}(B_r(S_{\ep,r}^{m-5}(u)))\leq Cr^5,
	\]
	which gives the second estimate of \eqref{emhm}. Moreover, note that $\dim_{\mathcal{H}}(S_{\ep,r}^{m-5})\le m-5$ implies that $\dim_{\mathcal{H}}({\rm sing}(u))\le m-5$. This gives a new proof for the estimate $\dim_{\mathcal{H}}({\rm sing}(u))\leq m-5$ without using the famous Federer dimension reduction method.
	
	 To prove the remaining estimate, simply observe that by Definition \ref{pppv} there holds
	\begin{align*}
	\{x\in B_1(0):r|\nabla u|+r^2|\nabla^2 u|+r^3|\nabla^3 u|+r^4|\nabla^4u|>1\}\su\{x\in B_1(0):r_u(x)<r\}.
	\end{align*}
	The proof is thus complete. 
\end{proof}

\begin{proof}[Proof Theorem \ref{thm85}]
It is similar to that of Theorem \ref{thm84} and hence is omitted.
\end{proof}

\appendix

\section{Classical stratification of singular set}\label{sec:appendix}
In this appendix, we present the classical stratification of singular set for biharmonic maps.

\begin{lemma}\label{inl}
	Under the assumption of Proposition \ref{prop: tangent measur}, there holds
	$$\lim_{j\ri\oo}H_{u_j}(y,r)=H_{v}(y,r),$$
	for all  $y\in \R^m$ and almost every $r>0$.
\end{lemma}
\begin{proof} 
	Fix any $0<r_0<R_0<\infty$ and write $A_{r_0,R_0}(y)=B_{R_0}(y)\backslash B_{r_0}(y)$. We first consider the case $y=0$ and aim to prove the following strong convergence:
	\[
	\pa_{X}|Du_{j}|^{2}+4|Du_{j}|^{2}-4\frac{|\pa_{X}u_{j}|^{2}}{|X|^{2}}\to\pa_{X}|Dv|^{2}+4|Dv|^{2}-4\frac{|\pa_{X}v|^{2}}{|X|^{2}}
	\]
	in $L^{1}(A_{r_0,R_0})$, where $X=x$. Indeed, since $Du_{j}\to Dv$ strongly
	in $L_{\loc}^{2}(\R^{m})$, we have
	\[
	4|Du_{j}|^{2}-4\frac{|\pa_{X}u_{j}|^{2}}{|X|^{2}}\to4|Dv|^{2}-4\frac{|\pa_{X}v|^{2}}{|X|^{2}}\qquad \text{in } L^{1}(A_{r_0,R_0}).
	\]
	To see that \begin{equation}\label{eq: partial strong convergence}
		\pa_{X}|Du_{j}|^{2}\to\pa_{X}|Dv|^{2}\qquad\text{in }L_{\loc}^{1}(\R^{m}),
	\end{equation}
	one needs to use the convergence \eqref{q6cb-2},
	which implies that $\pa_{X}D u_{j}\to-D v$ in $L_{\loc}^{2}(\R^{m})$.
	Thus 
	\[
	\pa_{X}|Du_{j}|^{2}=2D u_{j}\cdot\pa_{X}D u_{j}\to-2|D v|^{2}\qquad\text{in }L_{\loc}^{1}(\R^{m}).
	\]
	On the other hand, since $v$ is radially invariant, $\pa_{X}v=X\cdot D v=0$,
	it follows that $\pa_{X}\pa_{\al}v=x^\beta\pa_{\al\beta}v=-\pa_{\al}v$ for all $1\le\al\le m$.
	Hence
	\[
	-2|D v|^{2}=2D v\cdot\pa_{X}D v=\pa_{X}|Dv|^{2}.
	\]
	This proves \eqref{eq: partial strong convergence}.  Since $r_0$ and $R_0$ are arbitrary, this implies the claim. 
	
	In general, if $y\neq 0$, we note that 
	\[ H_{u_j}(y,r)=H_u(a+r_jy, r_j r)=H_{\tilde{u}_j}(0,1),\] 
	where $\tilde{u}_j(x)=u(a+r_j y+r_jrx)$, which is a scaling and small translation of $u(a+r_jrx)$. Thus the same result follows from the translation continuity of $W^{2,2}$-norms.
	The proof is complete.
\end{proof}

It follows from the above proof and the radial invariance of $v$ that 
\begin{equation}\label{eq: observation of H_v}
	H_{v}(0,r)=2r^{3-m}\int_{\pa B_r} |Dv|^2\equiv 2\int_{\pa B_1} |D_T v|^2,\qquad \text{for almost all }\,r>0,
\end{equation}
where
$D_T v$ is the tangent derivative of $v$ along $\pa B_1$.

Let $\mu=|\De v|^2dx+\nu$ be a tangent measure of a stationary biharmonic map $u$ at $a\in \Om$. For almost every $\rho>0$ and $y\in\R^m$, set
$$\Phi_{\mu}(y,\rho):=\rho^{4-m}\mu(B_\rho(y))+H_v(y,\rho)=\Phi_v(y,\rho)+\rho^{4-m}\nu(B_\rho(y)).$$

\begin{lemma}\label{kkkl}
	Let $u\in W^{2,2}(\Omega,N)$ be a stationary biharmonic map. Then $$y\in \reg(u)\quad \Longleftrightarrow\quad \Theta_{u}(y)=0.$$
\end{lemma}
\begin{proof} ($\Rightarrow$) If $y\in \reg(u)$, then $u$ is $C^2$ in a neighborhood of $y$ and so $\Theta_u(y)=0$.
	
	($\Leftarrow$)  Assume $\Theta_{u}(y)=0$. Let $\mu=|\De v|^2dx+\nu$ be a tangent measure of  $u$ at $y\in \Om$.  Then there exists a sequence $r_j\ri0$ such that the rescaled map $u_j\equiv u_{y,r_j}\rightharpoonup v$ in $W^{2,2}_{\loc},$ and so
	\begin{align}\label{ffft}
		\lim_{j\ri\oo}\Phi_{u_j}(0,r)=\lim_{j\ri\oo}\Phi_{u}(y,rr_j)=\lim_{r\ri0}\Phi_{u}(y,r)=\Theta_{u}(y)=0.
	\end{align}
	On the other hand, since $\mu$ is a bounded Radon measure, there have at most countably many $r>0$ such that $\mu(\pa B_r)>0$. Hence for almost every $r>0$, the weak convergence $|\Delta u_j|^2dx\wto \mu$ gives
	\[
	\lim_{j\ri\oo}\Phi_{u_j}(0,r)=\lim_{j\ri\oo}r^{4-m}\int_{B_r}|\Delta u_j|^2dx+H_{u_j}(0,r)=r^{4-m}\mu(B_r)+H_v(0,r)=\Phi_\mu(0,r),
	\]
	where we used Lemma \ref{inl}. Since $v$ is radially invariance and $\nu$ is a cone measure, there holds
	\[\Phi_\mu(0,r)\equiv \Theta_v(0)+\Theta_\nu(0),\qquad \text{for almost every }\,r>0.\] 
	Therefore, we find that
	\[\Theta_{u}(y)=\Theta_\mu(0,r)=\Theta_v(0)+\Theta_\nu(0)\]
	Moreover, due to the observation \eqref{eq: observation of H_v}, we have $\Theta_v(0)\ge 0$. Hence 
	\[\Theta_{u}(y)=0\quad\Longleftrightarrow\quad \Theta_v(0)=\Theta_\nu(0)=0.\]
	This implies that $v\equiv \text{constant}$. So
	$$r_i^{2-m}\int_{B_{r_i}(y)}|\na u|^2=\int_{B_1(0)}|\na u_{i}|^2\xrightarrow{i\to \infty} \int_{B_1(0)}|\na v|^2=0$$ 
	and by Lemma \ref{inl}
	\[H_{u_i}(0,r)\xrightarrow{i\to \infty} H_v(0,r)=0.\]
	By Lemma \ref{inl}, \eqref{ffft} and the above limit we have
	\begin{align*}
		\lim_{i\ri\oo}r_i^{4-m}\int_{B_{r_i}(y)}|\De u|^2=\lim_{i\ri\oo}\Phi_{u_{y,r_i}}(0,1)-\lim_{i\ri\oo}H_u(y,r_i)=0
	\end{align*}
	It follows from \cite[Lemma A.1]{Scheven-2008-ACV} that
	\begin{align*}
		r_i^{4-m}\int_{B_{r_i}(y)}|\na^2 u|^2\leq C_1r_i^{4-m}\int_{B_{2r_i}(y)}|\De u|^2+C_2 r_i^{2-m}\int_{B_{2r_i}(y)}|\na u|^2\xrightarrow{i\to\infty}0.
	\end{align*}
	Then Theorem \ref{thm3} gives that $y\in \reg(u)$ and the proof is complete.
\end{proof}



\begin{lemma}\label{ju01}
	Suppose $u\in W^{2,2}(\Omega,N)$ is a stationary biharmonic map and $v$ is a tangent map of $u$ at $a\in\Om.$ Then, for each $y\in \R^m$,  there exists a sequence $s_j\to \infty$ such that
	$$ \lim_{j\to \infty}H_v(y,s_j)=  H_v(0,1).$$
	The sequence ${s_j}$ can be taken uniformly for all $y\in B_R(0)$ for any fixed $R>0$. 
\end{lemma}
\begin{proof} Write $v_{y,r}=v(y+r\cdot)$. Then $H_{v}(y,r)=H_{v_{y,r}}(0,1)$.
	The 0-homogeneity of $v$ with respect to the origin implies that $v_{y,r}=v_{y/r,1}$.
	Since $v\in W^{2,2}_{\loc}(\R^m,N)$, we know from the translation continuity that
	\[
	\| v_{y,r}-v \|_{W^{2,2}(B_1\backslash B_{1/2})}=\| v_{y/r,1}-v \|_{W^{2,2}(B_1\backslash B_{1/2})}\xrightarrow{r\to \infty} 0. 
	\]
	This convergence is uniform for all  $y\in B_R(0)$ for any fixed  $R>0$. 
	In particular, for any $r\gg 1$, there exists $\lambda_r\in (1/2,1)$ such that
	\[
	H_{v_{y, r\lambda_r}}(0,1)-H_{v_{0,\lambda_r}}(0,1)=H_{v_{y, r}}(0,\lambda_r)-H_v(0,\lambda_r)\xrightarrow{r\to \infty} 0. 
	\]
	Since  $H_v(0,r)=H_v(0,1)$ for each $r>0$ by the 0-homogeneity of $v$, the conclusion follows by taking  $s_j =j \la_j.$
\end{proof}

Now we show that tangent measures also have monotonicity property as that of stationary biharmonic maps.
\begin{lemma}[Monotonicity of tangent measures]\label{kun}
	Suppose $u\in W^{2,2}(\Omega,N)$ is a stationary biharmonic map, $v$ is a tangent map of $u$ at a singular point $a\in {\rm sing}(u)$ and  $\mu=(v,\nu)$ is an  associated tangent measure as in Proposition \ref{prop: tangent measur}. 
	Then for every $y\in \R^m$, the function $\rho\mapsto \Phi_{\mu}(y,\rho)$ is nondecreasing and
	$$\Theta_\mu(y):=\lim_{\rho\ri0}\Phi_{\mu}(y,\rho ) \quad\text{exists and}\quad \Theta_\mu(y)\leq\Theta_\mu(0).$$ 
\end{lemma}
\begin{proof}
	By the definition of tangent map, there exists a sequence $r_j\ri 0$ (as $j\ri\oo$) such that
	$u_j\equiv u_{a,r_j}\wto v$ in $W^{2,2}_{\loc}(\R^m,N).$ By the monotonicity formula, we have
	\begin{align}\label{sinp}
		\Phi_{u_{j}}(y,r)-\Phi_{u_{j}}(y,\rho)\geq4(m-2)\int_{B_r(y)\backslash B_\rho(y)}\frac{|\pa_Xu_{j}(x)|^2}{|x-y|^m}dx
	\end{align}
	for each $y\in\Om$ and almost every $r>\rho>0,$ where $\pa_Xu_j(x)=(x-y)\cdot\na u_j(x).$
	On the other hand, by the definition of $\mu$ and Lemma \ref{inl}, we easily deduce that
	\begin{align*}
		\lim_{j\to\oo}\Phi_{u_{j}}(y,r)=\Phi_{\mu}(y,r).
	\end{align*}
	Taking $j\ri\oo$ in \eqref{sinp} we obtain
	\begin{align}\label{hhk}
		\Phi_{\mu}(y,r)-\Phi_{\mu}(y,\rho)\geq4(m-2)\int_{B_r(y)\backslash B_\rho(y)}\frac{|\pa_Xv(x)|^2}{|x-y|^m}dx\geq0,
	\end{align}
	Thus $\Phi_{\mu}(y,r)$ is nondecreasing in $r$  and $\Theta_{\mu}(y):=\lim_{\rho\ri0}\Phi_{\mu}(y,\rho)$ exists.
	
	Finally, applying Proposition \ref{prop: tangent measur} (iii) and Lemma \ref{ju01} yields
	\begin{align}\label{uukk}
		\Theta_{\mu}(y)&\leq\Phi_{\mu}(y,s_j) ={s_j}^{4-m}\mu(B_{s_j}(y))+H_v(y,s_j)\non\\
		&\leq\Big(\frac{s_j+|y|}{s_j}\Big)^{m-4}(s_j+|y|)^{4-m}\mu(B_{s_j+|y|}(0))+H_v(y,s_j)\non\\
		&\leq \Big(\frac{s_j+|y|}{s_j}\Big)^{m-4}\Big(\mu(B_1(0))+H_v(y,s_j)\Big)\\
		&\xrightarrow{j\to \infty} \mu(B_1(0))+H_v(0,1)=\Theta_{\mu}(0),\non
	\end{align}
	where $s_j$ is the sequence given by  Lemma \ref{ju01}.
\end{proof}

Now we are able to prove the classical stratification for singular set of $u$.
\begin{lemma}\label{ou}
	Suppose $a\in\sing(u) $ and $\mu$, $v$ and $\Phi_{\mu}$ are given as in Lemma {\rm\ref{kun}}. Define a set
	$$S(\mu)=\{y\in\R^m:\Theta_{\mu}(y)=\Theta_{\mu}(0)\}.$$
	Then $S(\mu)$ is a linear subspace of $\R^m$ and $\Theta_{\mu}$ is translation invariant along $S(\mu)$, i.e.
	$$\Theta_{\mu}(x+y)=\Theta_{\mu}(x)\qquad \text{for all }\,x\in\R^m,\,y\in S(\mu).$$ Moreover,  $S(\mu)\subset \sing(v)\cup \spt(\nu)$.
\end{lemma}
\begin{proof}
	For each $y\in \R^n$, we infer from \eqref{hhk} and \eqref{uukk} that
	\begin{equation}\label{uukk0}
		\begin{aligned} & 4(m-2)\int_{B_{s_{j}}(y)\backslash B_{\rho}(y)}\frac{\left|\partial_{X}v\right|^{2}}{|x-y|^{m}}dx+\Phi_{\mu}(y,\rho)\\
			& \le\Phi_{\mu}(y,s_{j})\le\Big(\frac{s_{j}+|y|}{s_{j}}\Big)^{m-4}\Big(\mu(B_{1}(0))
			+H_{v}(y,s_{j})\Big),
		\end{aligned}
	\end{equation}
	where $s_j$ is the sequence defined as in  Lemma \ref{ju01}. Then for any fixed $R>\rho$, we have
	\begin{align}\label{uukk1}
		4(m-2)\int_{B_R(y)\backslash B_{\rho}(y)}\frac{|\pa_Xv|^2}{|x-y|^m}dx+\Theta_{\mu}(y)&\leq\Theta_{\mu}(0).
	\end{align}
	Thus $y\in S(\mu)$ implies that $(x-y)\cdot\na v(x)=0$, i.e. $v(y+\lambda x)=v(y+x)$ for all $x\in B_R(y)\backslash B_{\rho}(y)$ and $\lambda>0$. Since $R>\rho$ are arbitrary, the previous equation holds for all $x\in \R^m$ and so $v$ is homogeneous with respect to both $0$ and $y$. For each $t\in \R$, choose $\lambda>0$ so that $t=\lambda-\lambda^{-1}$. Then
	\begin{align*}
		v(x)&=v(\lambda x)=v(y+\lambda x-y)=v(y+\lambda^{-2}(\lambda x-y))\\
		&=v(\lambda(y+\lambda^{-2}(\lambda x-y)))=v(x+ty),
	\end{align*}
	This shows that $v$ is translation invariant along the subspace $\R y$.
	
	Return to \eqref{uukk0}. Sending $s_j \to \infty$ yields
	\begin{align}\label{qqkk0}
		\Phi_{\mu}(y,\rho)\leq\Theta_{\mu}(0),\qquad \text{for all }\,\rho>0.
	\end{align}
	By the assumption $y\in S(\mu)$ and the monotonicity of $\Phi_{\mu}$, we deduce that
	\begin{align}\label{qq00}
		\Theta_\mu(0)=\Theta_\mu(y)\leq \Phi_{\mu}(y,\rho)\leq \Theta_{\mu}(0),\qquad \text{for all }\,\rho>0.
	\end{align}
	As a result, the homogeneity of $v$ with respect to $y$ implies that $\mu$ is scaling  invariant with respect to balls centered $y$. 
	
	As a result, we can further derive that $\mu$ satisfies $\mu(y+\la A)=\la^{m-4}\mu(A)$ for any $\la>0$ and Borel measurable set $A$, since $\mu$ is a cone measure. Then we have
	\begin{align*}
		\mu(B_1(x))&=\la^{4-m}\mu(y+B_\la(x)-y)\\
		&=\la^{(-2(4-m))}\la^{4-m}\mu(y+\la^{-2}(B_\la(x)-y))\\
		&=\la^{4-m}\la^{(-2(4-m))}\la^{4-m}\mu(\la(y+\la^{-2}(B_\la(x)-y)))\\
		&=\mu(B_1(x)+ty)=\mu(B_1(x+ty))=r^{4-m}\mu(B_r(x+ty)),
	\end{align*}
	for any $\la>0,x\in\R^m$ and $y\in S(\mu),$
	where $t=\lambda-\lambda^{-1}\in\R.$ This shows
	\begin{align}\label{qii}
		\Phi_{\mu}(x,r)=\Phi_{\mu}(x+ty,r).
	\end{align}
	So for any $a,b\in\R$ and $z_1,z_2\in S(\mu)$, there holds
	$$\Phi_{\mu}(x+az_1+bz_2,r)=\Phi_{\mu}(x+az_1,r)=\Phi_{\mu}(x,r).$$
	Letting $x=0$ and $r\ri0$ shows that $\Theta_{\mu}(az_1+bz_2)=\Theta_{\mu}(0).$ By the definition of $S(\mu)$, we have $az_1+bz_2\in S(\mu)$, which implies that $S(\mu)$ is a linear subspace of $\R^m$.
	Letting $t=1$ and $r\ri0$ in \eqref{qii}, we obtain $\Theta_{\mu}(x+y)=\Theta_{\mu}(x).$
	
	In order to show $S(\mu)\su \sing(v)\cup \spt(\nu)$,  we suppose now $y\in S(\mu)$  but $y\not\in \sing(v)\cup \spt(\nu)$. Then
	$$\Theta_\mu(y)=\lim_{r\to 0}\Phi_\mu(y,r)=0.$$
	This implies that $\Theta_u(a)=\Theta_\mu(0)=\Theta_\mu(y)=0$. So $a\in {\rm reg}(u)$ by Lemma \ref{kkkl}. Contradicting with the assumption that $a\in \sing(u)$. The proof is complete.
\end{proof}

\begin{corollary}\label{q2121}
	Suppose that $u\in W^{2,2}(\Omega,N)$ is a stationary biharmonic map, $a\in {\sing}(u)$ and $v$ is a tangent map of $u$ at $a$, and $S(\mu)$ is given by Lemma {\rm \ref{kun}}. Then
	$$\dim S(\mu)=m\Longleftrightarrow S(\mu)=\R^m\Longleftrightarrow \Phi_{\mu}\ \mbox{is a constant function }.$$
\end{corollary}
\begin{proof}
	The first equivalent relation is clear. For the second equivalence, note that by \eqref{qq00}, we have
	$$\dim S(\mu)=m\Longleftrightarrow \forall x,y\in\R^m,\Theta_\mu(x+y)=\Theta_\mu(x)\Longleftrightarrow\ \Phi_{\mu}\ \mbox{is a constant function }.$$
\end{proof}

From the theory of defect measure, we know that ${\rm dim} S(\mu)\leq m-4.$
No we define, for $j=0,1,\cdots,m-4$, the stratification of ${\rm sing}(u)$ by
$$\Si^j(u):=\{a\in{\rm sing}(u):\mbox{dim}(S(\mu))\leq j\ \mbox{ for all tangent measures } \mu  \mbox{ of } u  \mbox{ at } a  \}.$$
It is easy to see that
$$\Si^0\su\Si^1\su\cdots\su\Si^{m-4}=\Si^{m-3}=\Si^{m-2}=\Si^{m-1}={\rm sing}(u).$$
This is the so-called classical stratification of singular set of $u$. Note that if $u$ is minimizing biharmonic, then
\[
\Sigma^k(u)=S^k(u):=\{x\in \Omega:\text{no tangent maps of } u \text{ at } x \mbox{ is } (k+1)\mbox{-symmetric}\}.
\]

\end{document}